\documentclass[compmod]{rsl}
\gridframe{N}
\usepackage{amsmath,amssymb,amsxtra,latexsym,natbib,mathrsfs,url}
\usepackage[safe]{tipa}

\volume{0}
\issue{0}

\pyear{2010}
\pmonth{Month}
\doinu{10.1017/S1755020300000000}
\title[Review of Symbolic Logic]
{On the Syntax of Logic and Set Theory}
\author[L.T. Schoenbaum]
{LUCIUS T. SCHOENBAUM}
\affil{Louisiana State University\footnote{Received: August 25, 2008}}
\leftrunninghead{Lucius~T. Schoenbaum}
\rightrunninghead{On~the~Syntax~of~Logic~and~Set~Theory}

\newenvironment{axiom*}{\medskip\sc Axiom. \begin{it}}{\end{it}} 
\newenvironment{stantab}{\begin{tabbing}
\hspace*{.5in}\=\hspace*{.5in}\=\hspace*{.5in}\=\hspace*{.5in}\=\hspace*{.5in}\=\hspace*{.5in}\=
\hspace*{.5in}\=\hspace*{.5in}\=\hspace*{.5in}\=\hspace*{.5in}\=\hspace*{.5in}\=\hspace*{.5in}\= 
\kill}{\end{tabbing}} 

\newcommand{\lies}{\subset} 
\newcommand{\rise}{\supset} 
\newcommand{\lays}{\varsubsetneq} 
\newcommand{\rais}{\varsupsetneq} 
\newcommand{\tn}{\rightarrow} 
\newcommand{\lrn}{\leftarrow} 
\newcommand{\crn}{\leftrightarrow} 
\newcommand{\ein}{\ni} 
\newcommand{\nin}{\notin} 
\newcommand{\nein}{\not\ni} 

\newcommand{\formation}{\longmapsto} 
\newcommand{\nll}{\text{\upshape{\O}}}
\newcommand{\vstop}{{}'\mspace{1mu}} 
\newcommand{\istop}{\spdot} 
\newcommand{\cstop}{\mspace{2mu}\hat{}} 
\newcommand{\sbst}[3]{#1 [ #3 / #2 ]} 
\newcommand{\vas}[1]{\underline{#1}} 
\newcommand{\dn}{\bot} 
\newcommand{\kai}{\mspace{-1mu}\mathbin{{}_{\scriptscriptstyle{\setminus}}}\mspace{-3mu}} 
\DeclareMathOperator{\Sing}{\mathfrak{s}} 
\newcommand{\clas}[2]{\{{}^{#1}\, #2\}} 
\newcommand{\diff}{\setminus} 
\newcommand{\scomp}{\mathop{{}^{\scriptscriptstyle{\sim}}\mspace{-3mu}}}  
\newcommand{\dstnct}{\mspace{-2mu}\mathbin{{}_{\scriptscriptstyle{\Delta}}}\mspace{-2mu}} 
\newcommand{\ett}{\mspace{-5mu}\mathbin{{}_{\scriptstyle{\And}}}\mspace{-6mu}} 
\newcommand{\kain}[1]{\mathbin{{}_{{\scriptscriptstyle{\setminus}}_{#1}}}} 
\newcommand{\Singn}[1]{\Sing_{\scriptscriptstyle{#1}}} 
\newcommand{\inn}[1]{\in_{\scriptscriptstyle{#1}}} 
\newcommand{\liesn}[1]{\lies_{\scriptscriptstyle{#1}}} 
\newcommand{\midn}[1]{\mid_{\scriptscriptstyle{#1}}} 
\newcommand{\clasn}[3]{{}_{\scriptscriptstyle{#1}}\clas{#2}{#3}} 
\newcommand{\capn}[1]{\cap_{\scriptscriptstyle{#1}}} 
\newcommand{\SingI}{\Sing_{\scriptscriptstyle{1}}} 
\newcommand{\inI}{\inn{1}} 
\newcommand{\liesI}{\liesn{1}} 
\newcommand{\midI}{\midn{1}} 
\newcommand{\capI}{\capn{1}} 
\newcommand{\shl}{{}^{\text{\textendash}\mspace{2mu}}} 
\newcommand{\ashl}{{}^{\sim\mspace{2mu}}} 
\newcommand{\mdln}[1]{\mbox{\raisebox{1.2pt}{\u{}}}\mspace{-2.3mu}\overline{#1}} 
\newcommand{\mdlo}{\mdln{\phantom{Z}\;}}
\newcommand{\ssplus}{\:+\:} 
\newcommand{\ssmdln}[1]{\mbox{\raisebox{-0.52pt}{\u{}}}\mspace{-2.42mu}\overline{#1}}
\newcommand{\ssmdlo}{\ssmdln{\phantom{Z}\;}}
\DeclareMathOperator{\weight}{Wt} 
\DeclareMathOperator{\reduce}{red}
\DeclareMathOperator{\asg}{asg}
\DeclareMathOperator{\dec}{dec}
\DeclareMathOperator{\ev}{ev}
\newcommand{\call}[1]{\mathop{\forall #1 \,}} 
\newcommand{\calln}[2]{\mathop{\forall_{\scriptscriptstyle{#1}} #2 \;}} 
\newcommand{\there}[1]{\mathop{\exists #1 \,}} 
\newcommand{\theren}[2]{\mathop{\exists_{\scriptscriptstyle{#1}} #2 \;}} 
\newcommand{\therenll}[1]{\mathop{\exists^{\scriptscriptstyle{\nll}} #1 \;}} 
\newcommand{\Z}{\mathbb{Z}} 
\newcommand{\ksi}{\xi}	

\begin{document}
\maketitle

\begin{abstract}
We introduce an extension of the propositional calculus to include abstracts of predicates and quantifiers, employing a single rule along with a novel comprehension schema and a principle of extensionality, which are substituted for the Bernays postulates for quantifiers and the comprehension schemata of ZF and other set theories.  We prove that it is consistent in any finite Boolean subset lattice.  We investigate the antinomies of Russell, Cantor, Burali-Forti, and others, and discuss the relationship of the system to other set theoretic systems ZF, NBG, and NF.  We discuss two methods of axiomatizing higher order quantification and abstraction, and then very briefly discuss the application of one of these methods to areas of mathematics outside of logic.  
\end{abstract}


\section{Introduction}\label{S:1}

This work is an extended investigation of the idea of regarding the set $\{a,b\}$ not as the extension of a class (in this case, the class whose existence is guaranteed in ZF by the Pairing axiom), but as a mathematical product $a \star b$, an approach which, as we will see, naturally gives rise to a hierarchy of products $\star_1, \star_2, \star_3 \ldots$ in some ways comparable to the recursive hierarchy of arithmetical operations $+,\cdot,$ etc.  A major goal is to obtain a well-defined notion of the abstract $\{ x \mid \phi(x) \}$ of the predicate $\phi$ such that, informally stated, one has $\{ x \mid \phi(x) \lor \psi(x) \} = \{ x \mid \phi(x) \} \star \{ x \mid \psi(x) \}$, and other theorems of the same kind.  There are two facts to underline:
\begin{enumerate}[(ii).]
\renewcommand{\theenumi}{(\roman{enumi})}
	\item In order to form such abstract classes, one is forced to introduce some kind of notion of $\star$-primality.  Havoc would result if $\star$-divisible objects could be ``included" or ``added" to the class by satisfying the condition $\phi$, because this property might not necessarily be inherited by its $\star$-divisors.
	\item In order to form such abstract classes, the empty set must acquire a status unlike that of ordinary sets.  Were such a class to be empty, $\{x \mid \phi(x) \} = \nll$, it must be well defined whether it is meant that $\phi(\nll), \lnot\phi(\nll)$, or perhaps neither.  
\end{enumerate}
It is not surprising that difficulties arise once, so to speak, the animals have all been released from their cages---once, that is, braces are eliminated and totalities give way to trickier (and more fluid) multiplicities.  What is remarkable is that it seems that surpassing these obstacles, as a straightforward consequence, yields a naive comprehension principle, since two critical paradoxes, those of Russell and of Cantor, are correlated with the two points above.  Cantor's paradox, which arises due to the unwieldy growth of power sets, can be avoided if the power set is not allowed to be $\star$-prime.  Russell's paradox is avoided if we do not grant $\star$-primality to the empty set.  \pagebreak 

There is a second idea, closely related to the first, involved in our investigation.  By building set theory not as the semantics of first order logic, but (through a kind of ``radical Henkinization") syntactically of a piece with it, we demonstrate that the basic philosophical principles of formalism can be applied in a setting further removed from logicist foundations than that setting established by Hilbert and his circle during the formative years of mathematical logic.  Through a slight modification of the model-theoretic boundary between syntactic and semantic spaces, we invest formulas with location among model-theoretic individuals, allowing a number of simplifying reductions of the main artifacts of logic to be implemented.  This appears to carry interesting effects through to more advanced levels of analysis in proof theory and model theory (\S 3), set theory (\S 4), and in other areas of mathematics (\S 5).  From a broader vantage, it gives rise to an organization or \emph{Weltanschauung} for mathematics which may be of interest to scientists and philosophers in areas related to mathematical foundations.  

Our presentation is essentially self-contained, in order to address a broad audience inclusive of researchers (philosophers and mathematicians alike) in areas closely related to logic who may not have a background in logic comparable to those of specialists in the field.  We would be most content if our work were considered interdisciplinary, and sincerely hope that our presentation contains something of interest to much-admired target readers in a plurality of disciplines.

\section{Elements}\label{S:2}

\subsection{Preliminaries}\label{S:2.1}
Our treatment of propositional calculus in this section will be brief.  We study the domain of \emph{objects}, or (synonymously) \emph{expressions}.  Anything that ``exists", we say, is an object---for example, propositions, desks, and numbers are objects.  A handful of objects are designated as \emph{signs}, these are:
\begin{equation*}
 \nll\quad 	\vstop\quad 	\istop\quad 	\cstop\quad 	\tn\quad	\kai\quad 	\cap\quad
\end{equation*}
The symbol $\nll$ is the system's basic \emph{constant}.  The next three symbols of the list generate new objects.  An object of the form $\alpha\vstop$ for some object $\alpha$ is a \emph{variable}, an object of the form $\alpha\istop$ is an \emph{index}, and a \emph{constant} is, in addition to $\nll$, any object of the form $\alpha\cstop$.  We call the first handful of variables $\nll\vstop, \nll\vstop\vstop, \nll\vstop\vstop\vstop,$ etc. with the standard abbreviations $a, b, c, x, y, z, A, B, C, $ etc.; this assignment is taken to be precise, but need not be specified here.  In our presentation, when we have a need for metamathematical variables, we use Greek letters $\alpha, \beta, \gamma,$ etc.  We call the three unary connectives $\vstop, \istop, \cstop$  \emph{stops} or \emph{logical stops}.  This name is chosen since they are meant to ``stop" the progress of a logical substitution (defined below).  The index stop can be put aside for now; we will not see it again until \S 3.  

The symbol $\tn$ is the fundamental \emph{relation}, referred to as \emph{containment} or \emph{universal containment}.  The symbols $\kai$ and $\cap$ are (binary) \emph{operations}.  Table \ref{F:frules} displays the system's basic formation rules.  These rules allow, for any objects $\alpha$, $\beta$, $\gamma$, the formation of a new object, having $\alpha$, $\beta$, $\gamma$ as \emph{subobjects} or \emph{subexpressions}; this language is heritable in the obvious sense.  The object $(\alpha \kai \beta)$ is the \emph{conjunction}, \emph{union}, or \emph{federation} of $\alpha$ and $\beta$.\footnote{We have searched for an apt name and symbol for the product ($\alpha \kai \beta$): a universal expression of objects-in-multiplicity, akin to ordinary addition but idempotent and, in practice, more basic.  The symbol employed is a formal comma modeled after the commas in the notation $\{a,b\}$ for the pair and the Gentzen-Kleene-Rosser notation $\Gamma, A \vdash B$ for the carriage of assumptions.  The use of an ordinary comma is tempting but troublesome, so one considers using one or the other of $\cup$ and $\land$.  However, because the product should carry neither logical nor set-theoretical connotation in every instance, stultifying and plainly odd-looking expressions arise as a result.  The symbol $(\kai\,)$  smooths over these distinctions, and thus helps us achieve the desired conceptual view.  Our approach is to use the terminology \emph{federation}, \emph{union}, and \emph{conjunction} (and three distinct symbols) to refer to instances of the same operation appearing at three distinct positions in the order of operations (see footnote, p. \pageref{OoOp}).  It is hoped that a new symbol and a new term will aid the reader and avoid confusion to the greatest possible extent.}  The object $(\alpha \cap \beta)$ is the \emph{intersection} of $\alpha$ and $\beta$.  We shall suppress these parentheses when they are unnecessary.  

We call any object generated from $\nll$ using only the given formation rules (and those that may be added later) \emph{formal objects}, or \emph{well-formed} objects.  Objects which are not well-formed are \emph{informal objects}, and include chairs, desks, letters of the alphabet (when not intended as abbreviations of formal objects), and strings of signs which are not well-formed.  In this paper, we shall loosely allow consideration of informal or ``ordinary" objects as objects of the system, while giving precedence to the view in which all but the formal objects are laid aside or ignored.  Thus the terms \emph{object} and \emph{expression} will hereafter refer to formal objects.  

\begin{table}[t]\caption{Formation Rules (objects $\alpha$, $\beta$)}\label{F:frules}
\begin{center}
\begin{oldtabular}{c}
\begin{oldtabular}{l}
$\nll \formation \nll$ \\
$\hspace{0.5pt}\alpha\hspace{0.6pt} \formation (\alpha)\cstop$ \\
$\hspace{0.5pt}\alpha\hspace{0.6pt} \formation (\alpha)\vstop$ \\
\end{oldtabular} \qquad
\begin{oldtabular}{l}
$\alpha \, \beta \formation (\alpha \kai \beta)$ \\
$\alpha \, \beta \formation (\alpha \tn \beta)$ \\
$\alpha \, \beta \formation (\alpha \cap \beta)$ \\
\end{oldtabular} \\
\end{oldtabular}
\end{center}
\end{table}

We make the following definitions.  For $\nll\cstop$ we write $\dn$.  For $a \tn \dn$ we write $\lnot a$, the \emph{negation} of $a$.  For $\nll \tn a$ we write $\vas{a}$, the \emph{trivialization} of $a$.\footnote{This notation mirrors the notation $\bar{a}$ for the negation of $a$ in, e.g., \citet{Hilbert1928}.}  We write $a \lor b$ for $\vas{\vas{a}\cap\vas{b}}$, the \emph{disjunction} of $a$ and $b$.  This expression is classically equivalent to the expression $\lnot(\lnot a \kai \lnot b)$; however, the former expression remains serviceable in the intuitionistic case.\footnote{Thus $\cap$ can replace Heyting's connective $\lor$.}  For $b \tn a$ we may write $a \lrn b$, and for $(a \tn b) \kai (a \lrn b)$ we write $a = b$.  

We call expressions which may be put in the form $a \tn b$\footnote{Given an object $\ksi$, that is, we may find objects $\alpha$ and $\beta$ such that $\ksi = (\alpha \tn \beta)$.} \emph{formulas}; an object that is not a formula is said to be \emph{concrete}. Expressions of the form $a = b$ may also be called \emph{equations} or \emph{equalities} in addition to being called formulas, as usual.  In order to improve readability and respect normal usage, we will usually use the synonymous symbol $\crn$ instead of $=$ when we denote an equality between formulas.  We may also similarly write $a \land b$ in place of $a \kai b$ in the case when $a$ and $b$ are both formulas.  In its role, this technique is comparable to the practice (common, for example, in physics) of using modified parentheses (brackets, braces) in lengthy expressions to highlight syntactic units.  It is also comparable in its role to the system of dots devised by Peano and featured in the \emph{Principia Mathematica}.  Here, because our work is not too involved, and in order to gain familiarity with the machine language (so to speak), we shall use $\land$ sparingly.  

We must pause to clarify the aforementioned notion of logical substitution.   First, an instance of an object $\beta$ in the syntax of an object $\alpha$ may be one of two kinds: if $\beta$ is, in this instance, a subobject of $\alpha$, it is said to \emph{appear} in $\alpha$, and otherwise to merely \emph{occur} in $\alpha$.  We denote the object $\alpha$ modified so that every appearance (note) of the variable $x$ is replaced by the object $\beta$ by $\sbst{\alpha}{x}{\beta}$.  We refer to $x$ as the \emph{target variable}, and to $\beta$ as the \emph{substituend} of the substitution.  The reader can verify that this is a well-defined object for its inputs $\alpha$, $\beta$, and $x$, and that no expression which occurs within (``is guarded by") an index, variable, or constant stop is modified by any substitution.  If there is no risk of confusion, in place of $\sbst{\alpha}{x}{\beta}$, we may sometimes write simply $\alpha(\beta)$.  

\begin{table}[t]\caption{Postulates of Section 2 (objects $\alpha$, $\beta$, variables $a, b, c, d, x$).}\label{F:rax1}
\begin{center}
\begin{oldtabular}{c}
		\begin{oldtabular}{cc}
		
		R1.  Exchange.  \qquad\qquad\qquad\quad & R2.  Substitution. \\
		$\begin{array}[b]{l}
			\,\alpha \\ 
			\,\alpha \tn \beta \\
			 \,\beta						\rule{0pt}{10pt}\\
			\raisebox{24.5pt}[0pt][0pt]{$\begin{array}{c}\phantom{xxxx}\\ \hline \end{array}$}
		\end{array}$ \qquad\qquad\qquad\quad
		&
		$\begin{array}[b]{l}
			\\
			\,\alpha \\
			\,\sbst{\alpha}{x}{\beta} 			\rule{0pt}{10pt}\\
			\raisebox{24.5pt}[0pt][0pt]{$\begin{array}{c}\phantom{xxx\;}\\ \hline \end{array}$}
		\end{array}$
		
		\end{oldtabular}\\

\begin{oldtabular}{l}
	A1.  $((c \tn a) \kai (c \tn (a \tn b))) \tn (c \tn b).$ \\
	A2.  $((d \tn (a \tn b)) \kai (d \tn (b \tn c))) \tn (d \tn (a \tn c)).$ \\
	A3.  $((c \tn a) \kai (c \tn b)) \tn (c \tn (a \kai b)).$ \\
	\\
	\begin{oldtabular}{l}
		A4a.  $(a \kai b) \tn a.$ \\
		A4b.  $(a \kai b) \tn b.$
	\end{oldtabular} 
	\qquad \qquad \qquad
	\begin{oldtabular}{l}
		A5a.  $a \tn (a \cap b).$ \\
		A5b.  $b \tn (a \cap b).$ 
	\end{oldtabular} \\
	\\
	A6.  $((a \tn c) \kai (b \tn c)) \tn ((a \cap b) \tn c).$ \\
\end{oldtabular}\\
\\
\begin{oldtabular}[t]{l}
	A7.  $a \tn a.$ \\
	A8.  $a \tn \nll.$ \\
	A9.  $\dn \tn a.$
\end{oldtabular}
\quad \qquad \qquad
\begin{oldtabular}[t]{l}
	A10.  $(\nll \tn a) \tn ((\nll \tn b) \tn (a \kai b)))$ \\
	A11.  $(a \tn b) \tn (\nll \tn (a \tn b)).$ \\
	A12.  $(((a \tn b)\tn \dn)\tn \dn) \tn (a \tn b).$
\end{oldtabular}\\		
\end{oldtabular}
\end{center}
\end{table}

A \emph{formal proof} $\Pi = (\alpha_1 , \ldots , \alpha_n)$ is a finite numbered list of substitutions, exchanges, and axioms.  These \emph{steps} of the proof are usually written in descending order down the page along with a citation or attribution column (as in the formal proof of Lemma 2.1 
below).  If it is possible for an expression to appear in a formal proof, we say that the expression may be \emph{obtained} or \emph{deduced}, and we call $(\alpha_1, \ldots, \alpha_n)$ a \emph{formal proof of $\alpha_n$}.  Rules R1-2 and axioms A1-12 are presented in Table \ref{F:rax1}.  They are styled \emph{postulates}; a postulate is a rule or an axiom, according to \citet{Kleene1952}.  They are based primarily on the systems of \citet{Heyting1930} (see \citet[pp. 311-27]{Mancosu1998}) and Kleene, \emph{op. cit.}, but see also \citet{Hilbert1927}, and \citet[pp. 55ff.]{Hilbert1928}.  The system is what is known in proof theory as a substitution Frege system (see, e.g., \citet{Pudlak1998}).  The rule of \emph{Combination}:
\begin{equation*}
\begin{array}{c}
	\alpha \\
	\beta	\\
	\alpha \kai \beta			\rule{0pt}{10pt}\\
	\raisebox{24.5pt}[0pt][0pt]{$\begin{array}{c}\phantom{xx\,}\\ \hline \end{array}$}
\end{array}
\end{equation*}
\noindent is immediately derived from A10 and A11, for formulas $\alpha$ and $\beta$.  

\subsection{Basics}\label{S:2.2}
Formal proofs are greatly simplified through the use of \emph{assumptions}.  An arbitrary object (normally a formula) $\alpha$ is taken or \emph{assumed} at a given step, and subsequent rules may combine relying on the assumption as if it were proven, until the assumption is \emph{discharged}.  At the discharge step the conclusion of the immediately preceding step is rewritten under the hypothesis of the assumption.  A proof with assumptions is not complete until all assumptions have been discharged.  That the use of assumptions does not allow any new formulas to be proven is a result known as the deduction theorem, which we will now prove.  We first require a lemma.  

\begin{lem}\label{T:FK}
$\vas{a} \tn (b \tn a).$
\end{lem}
\begin{proof*} 
\begin{stantab}
	$1.\;   a \tn \nll$						\>\>\>\>\>\>\> A8  \\
	$2.\;  (\nll \tn a) \tn \nll$				\>\>\>\>\>\>\> $\sbst{(1)}{a}{(\nll \tn a)}$ \\
	$3.\;   a \tn \nll$	 					\>\>\>\>\>\>\>A8 \\
	$4.\;   b \tn \nll$						\>\>\>\>\>\>\> $\sbst{(3)}{a}{b}$ \\
	$5.\;  (a \tn b) \tn (\nll \tn (a \tn b))$		\>\>\>\>\>\>\> A11 \\
	$6.\; (a \tn \nll) \tn (\nll \tn (a \tn \nll))$		\>\>\>\>\>\>\> $\sbst{(5)}{b}{\nll}$ \\
	$7.\; (b \tn \nll) \tn (\nll \tn (b \tn \nll))$  	\>\>\>\>\>\>\> $\sbst{(6)}{a}{b}$ \\
	$8.\; \nll \tn (b \tn \nll)$				\>\>\>\>\>\>\> Exch. 4, 7 \\
	$9.\; (a \tn b) \tn ((b \tn c) \tn (a \tn c))$		\>\>\>\>\>\>\> Lemma (A4a, A4b,\\
									\>\>\>\>\>\>\> Comb., A2, Exch.) \\
	$10.\; ((\nll \tn a) \tn b) \tn ((b \tn c) \tn ((\nll \tn a) \tn c))$ \>\>\>\>\>\>\> $\sbst{(9)}{a}{\nll \tn a}$ \\
	$11.\; ((\nll \tn a) \tn \nll) \tn ((\nll \tn c) \tn ((\nll \tn a) \tn c))$ \>\>\>\>\>\>\> $\sbst{(10)}{b}{\nll}$ \\
	$12.\; ((\nll \tn a) \tn \nll)$ \\
	\qquad $\tn ((\nll \tn (b \tn \nll)) \tn ((\nll \tn a) \tn (b \tn \nll)))$ \>\>\>\>\>\>\> $\sbst{(11)}{c}{(b \tn \nll)}$ \\
	$13.\; (\nll \tn (b \tn \nll)) \tn ((\nll \tn a) \tn (b \tn \nll))$\>\>\>\>\>\>\> Exch. 2, 12 \\
	$14.\; (\nll \tn a) \tn (b \tn \nll)$			\>\>\>\>\>\>\> Exch. 8, 13 \\
	$15.\; a \tn a$						\>\>\>\>\>\>\> A7 \\
	$16.\; (\nll \tn a) \tn (\nll \tn a)$			\>\>\>\>\>\>\> $\sbst{(15)}{a}{(\nll \tn a)}$ \\
	$17.\; ((d \tn (a \tn b) \kai (d \tn (b \tn c))) \tn (d \tn (a \tn c))$ \>\>\>\>\>\>\> A2 \\
	$18.\;$---							\>\>\>\>\>\>\> $\sbst{(17)}{c}{(\nll \tn a)}$ \\
	$19.\;$---							\>\>\>\>\>\>\> $\sbst{(18)}{b}{\nll}$ \\
	$20.\;$---							\>\>\>\>\>\>\> $\sbst{(19)}{a}{b}$ \\
	$21.\;$---							\>\>\>\>\>\>\> $\sbst{(20)}{c}{a}$ \\
	$22.\; ((\nll \tn a) \tn (b \tn \nll)) \kai ((\nll \tn a) \tn (\nll \tn a))$ \>\>\>\>\>\>\> Comb. 14, 16 \\
	$23.\; (\nll \tn a) \tn (b \tn a)	$			\>\>\>\>\>\>\> Exch. 22, 21
\end{stantab}
\end{proof*}

Because of A11, and because A1 implies that $(\nll \tn a) \tn a$, the equation
\begin{equation}\label{E:visaca}
\vas{a} = a
\end{equation}
is characteristic of formulas.  Thus, Lemma 2.1 
becomes, for all formulas $a$ and $b$, formula (1) of \citet{Frege1879}, 
\begin{equation}\label{E:FKpq}
a \tn (b \tn a).
\end{equation}
The connection between (\ref{E:FKpq}) and assumption mechanisms was noticed early on (see \citet[p. 465]{Heij}).  However, it does not obtain when $a$ and $b$ range over all the objects of our system, as can easily be checked by letting $a$ in (\ref{E:FKpq}) be concrete.\footnote{Provided that one has already crossed over to the intended conceptual view, this point is obvious.  We are considering logic and set theory as they might be conceived in a single common space or worldview in which a commitment exists to a general notion of objecthood encompassing both language and its referents.  Consider, then, the content of (\ref{E:FKpq}) set-theoretically, e.g., in the case when $a$ and $b$ are ordinary objects (desks, chairs, etc.) and the relation is one of ordinary containment.  In our system concrete objects live in the midrange or midst of the lattice.  Viewed logically, the lattice encompasses an infinite array of truth values; a concrete expression---a desk, a number, a geometric figure, etc.---is an object whose truth value is neither true, nor false, but rather a thing unique to the object itself.  As the formalism is developed in the sequel, the load on the symbol $\tn$ will be decreased by the introduction of the symbols $\rise, \lies, \in$ etc.---though the distinction between these symbols and the family $\tn, \lrn,$ etc. shall solely involve properties of their arguments.  The formal unity between the concepts of implication and containment, whose naturality we propose in this section and whose fundamental importance to the system should become clear as we proceed, rests undisturbed throughout (echoing in certain respects an early system of Quine, see footnote, p. \pageref{inclusionandabstraction}).}

\begin{thm}[Deduction Theorem]\label{T:deduct}
A proof with assumptions in which all variables appearing in assumptions are held constant, and in which all assumptions are assumed and discharged in last-in-first-out order, may be converted to a formal proof.
\end{thm}
\begin{proof*}
Let $(\alpha_1, \ldots, \alpha_n)$ be such a proof with assumptions.  Consider a sequence of steps $(a_i, \ldots, $ $a_{i+j})$, $1 \leq i < i+j \leq n$,  from an assumption step, $\alpha_i$, to the first subsequent discharge step, $\alpha_{i+j}$, inclusive, containing no intervening assumption steps.  Replacing the entire string with the new string $(\alpha_i \tn \alpha_{i+1}, \alpha_i \tn \alpha_{i+2}, \ldots, \alpha_{i} \tn \alpha_{i+j-1})$, it is clear that every prior use of the rule of Exchange can be carried out under the lingering hypothesis $\alpha_i$ by inserting the Axiom A1, instantiating it as desired,\footnote{As a minor additional hypothesis, one must let variables appearing in assumptions be distinct from those appearing in the axioms.} and applying Combination and Exchange.  Every Substitution can be executed under the assumption since, by hypothesis, all variables are held constant.  All of the axioms, finally, can themselves be written under the assumption of $\alpha_i$ by using Lemma 2.1 
and Exchange, and we obtain a new formal proof.  If it has assumptions, we iterate the process described.  If we continue until we have treated all the assumptions in turn, we obtain a formal proof of $\alpha_n$.  
\end{proof*}

The stipulation that assumed variables be held constant is perfunctory at this level; there is no motivation to vary assumptions until quantifiers (and/or classes) are introduced.  All naturally arising assumptions may in fact be discharged in any order, but one must exercise a bit of caution, since
\begin{equation}
a \tn (b \tn c) \quad \tn \quad (a \kai b) \tn c \quad \tn \quad \vas{a} \tn (b \tn c),
\end{equation}
and these expressions are not equivalent in general,\footnote{If your wallet is empty, I will give you my coat.  But I am not inclined, otherwise, to give you my coat for your wallet.} but are seen to be equivalent when $a$ is a formula.  A few details remain which we omit showing that we may generalize these principles, and therefore make and manipulate assumptions naturally, provided that they are formulas.  

The following relations easily obtain: conjunctions and intersections of formulas are formulas.  The converses of axioms 3, 6, 11, and 12.  Of union/conjunction/ federation $(\kai\,)$: associativity, commutativity, idempotence, identity $\nll$, zero element $\dn$, as well as the relations
\begin{align}
	(a \tn b) &\tn (a \kai c \tn b \kai c)	\quad		\text{(isotonicity)} \\
	(a = b) &\tn (a \kai c = b \kai c)	\quad		\text{(substitutivity)}
\end{align}
Of intersection $\cap$: associativity, commutativity, idempotence, identity $\dn$, zero $\nll$, and isotonicity (implying substitutivity\footnote{The term \emph{isotone} (and the dual \emph{antitone}) is found, e.g., in \citet{Birkhoff1940}; \emph{substitutivity} appears, e.g., in \citet{QuineMLrevised,Bernays1958}.  One may show that antecedents are antitone, that double antecedents are again isotone, etc.}).  The same laws therefore hold also for disjunction $\lor$.  Finally, both of the distributive laws hold, yielding a distributive lattice of objects under federation and intersection, and a Boolean ring of formulas (under conjunction and disjunction, or equivalently, federation/union and intersection) isomorphic to $\Z_2$.

\subsection{Remarks}\label{S:2.3}
All students of modern mathematics are well accustomed to the assignment of formulas and equations to some value or thing that is in some sense the storing ground of true formulas, and another which is the home and identifier of false formulas.  It is a truth which is all the more remarkable given the advanced technology, both industrial and intellectual, available to scientists and philosophers who investigate these instruments in our own age, that whether they are taken to be abstract forms, letters, or empirical entities, they have up to now resisted all efforts to make them rigorously understood and canonized.  We do not contest these efforts, which lie well beyond our ken.  We follow formalists and intuitionists in distinguishing mathematical discourse from ordinary linguistic discourse, and place before ourselves only the very limited task of calibrating and arranging adequate grounds for the former, laying aside questions surrounding the latter.  We make no claims here concerning the concepts of truth and falsehood as they arise in ordinary language.  

Table \ref{F:pcstl} illustrates the syntactic medium in which we are developing logic and set theory.   This arrangement may surprise  some readers, many of whom may never have paused before over the evocative relations which naturally arise in proof theory (such as $\nll \vdash \Delta$, denoting that the set of formulas $\Delta$ follows under \emph{no assumptions}).  The author readily agrees that the arrangement is conceptually surprising (it stirs the imagination even after long familiarity); however, having weighed the matter carefully, he finds arguments which speak for the system as it has been defined, and to condemn it, neither any falsifier, nor such a drastic departure from the norm in foundations.  Logic itself is surprising and elusive---it receives these qualities naturally and instantly from the concepts it traffics in.  Assigning those formulas which are formally true to the formal constant $\nll$ identifying empty and trivial structures is no more remarkable than assigning them to a letter $\top$, ensconced in a detached logical space.  The formulas, in either case, remain what they are; we encounter them, and they inspire us, in the same way.  The author believes that it is rare indeed to find an a priori basis for criticism of variations in formal orientation which is invulnerable to a challenge on a priori grounds. One must be practical and, as Hilbert implored us to do, measure theories by their fruits.  Given the inconclusive suggestion of patterns in the structure of formal and informal thought, we should be neglecting the due diligence of science to leave either possibility unexamined, in case one orientation should offer distinct utilitarian advantages over its alternative.  The author cautions the reader to carefully entertain a wide-ranging field of evidence before adopting the Einsteinian critique that the system is simple, yes---and \emph{too} simple.  So far, we have seen a number of ways in which the present system recommends itself: a compact set of basic symbols, accessible definitions, the felicitous distinction between the formula and the concrete term, as well as the convenience of a unified intuitive presentation.  What remains to be discovered is, of course, unknown to the author.  In the next section, the reader shall have the opportunity to consider the effects of the novel approach upon the fundamental devices of set theory and model theory.  

\begin{table}[t]\caption{Propositional calculus in a set-theoretical lattice.}\label{F:pcstl}
\begin{center}
\begin{oldtabular}{cccc}
\; Truth Table \quad	& Name	 	& 	Generalization	& Formula							\\
TTTT &	TRUE	&	empty/trivial set 			& $\nll$							\\
FTTT &	NAND	&	$\{a,b\}$ is a covering of $\dn$	& $\lnot(a \kai b)$ 					\\
TFTT &			&	containment				& $a \tn b$ 						\\
TTFT &			&	 containment 				& $b \tn a$						\\
TTTF &	OR		&	 at least one set is empty/trivial	& $\vas{\vas{a} \cap \vas{b}}$			\\
FFTT &			&	 $\{a\}$ covers $\dn$			& $\lnot a$							\\
FTFT &			&	 $\{b\}$ covers $\dn$			& $\lnot b$							\\
FTTF &	XOR		&	 disjointness 				& $\vas{a \cap b} \kai (\nll \neq a \kai b)$	\\
TFTF &			&	 $b$ itself					& $b$							\\
TTFF &			&	 $a$ itself					& $a$							\\
FFFT &	NOR		&	 \quad all sets are nonempty/nontrivial \qquad	& $\lnot\vas{\vas{a} \cap \vas{b}}$		\\
FFTF &			&	 strict containment 			& $b \tn a \kai \lnot(a \tn b)$		\\
FTFF &			&	 strict containment 			& $a \tn b \kai \lnot(b \tn a)$				\\
TFFF &	AND		&	 union 					& $a \kai b$						\\
FFFF	 & \quad\; FALSE \;\;\;\; &	 absurd/undefined 		& $\dn$							\\
\end{oldtabular}\\
\begin{oldtabular}{c}
\raisebox{180pt}[0pt][0pt]{$\begin{array}{c}\phantom{xxxxxxxxxxxxxxxxxxxxxxxxxxxxxxxxxxxxxxxxxxxxxxxxxxxxxxxxxxxx}\\ \hline \end{array}$}
\end{oldtabular}\vspace{-12pt}
\end{center}
\end{table}

Aside from these remarks (and with respectful apologies to those in the philosophical domain), the author does not wish to suggest here a single philosophical vantage from which the formal system might best be viewed.  According to strict formalism, the system is the article found on the page, dissociated from all semantic interpretations, and whatever metaphysics it might inspire is irrelevant to its use in practice.  The formalist, as Carnap once said, is ``among the most pacifistic of mankind," able to partake in the methods of all others, and the author has, in principle, no disagreements with him.  We believe that one could incorporate the system into intuitionism by introducing, to begin with, the symbol $\nll$ as the empty hand---the moment before the unfolding construction of proof---while introducing the symbol $\dn$ to represent an impossible object of reflection, or what is called in \citet{Godel1944} ``the notion of `something' in an unrestricted sense".  One may then proceed to define the system of Heyting.\footnote{The intuitionistic case is interesting in its details, but requires additional work; we intend to discuss it elsewhere.  The opposite assignment conflicts with intuitionism, for consider: truth is obviously constructive, while the domain in which construction occurs cannot be entirely so.}  With the early work of Wittgenstein, there is a shared goal of a comprehensive \emph{Weltanschauung} for analytic thought, and with sustained effort along these lines there might be developed a framework bearing out some of his views on the limitations of language and the world of coherent ideas.  The relationship of these ideas to Platonism and its many later manifestations is also quite striking.  The contradiction-as-unity ($\nll = \dn$) has some curious and remarkable parallels to concepts explored by several ancient Greek philosophers including Plato, later figures such as, for example, Plotinus, Cusanus, and Pascal, and modern figures such as Cantor, Brouwer, Whitehead, and many others.  We believe considerable opportunity remains to discuss the system in light of other leading schools.  The relationship with Cantor's thought is further touched upon in \S 4.2, 
where we shall see that the novel approach adds a note of prescience to one of his more rarely-cited ideas.   

\section{Predicate Calculus and Set Theoretical Abstraction}\label{S:3}

To first summarily recount our progress in \S 2: 
to an unformalized intuitive space we have added the presence of \emph{objects}, allowed that these objects may be collected into multiplicities (also considered objects), defined the intersections of such multiplicities (again, objects), and developed a system generating formulas (also objects) describing the containment relationship, if any, holding between two given objects.  To this world of things, we wish to apply set-theoretical abstraction.  This demands that we make a few incomplete remarks, before we begin, on the general question of how best to introduce sets of elements satisfying a stated condition---abstracts, or what most authors call \emph{classes}---to a system like ours.  

\subsection{Individuals}\label{S:3.1}

In an early attempt to develop the class rigorously for analytic purposes, \citet{Peano1889} refers to the formation of classes as \emph{inversion}, apparently on the notion that if a formula were taken as a mapping into the truth values 0 and 1, the class derived from the formula would be its kernel. {\em Prima facie}, this seems quite natural, perhaps even obvious, but it runs into difficulties having to do with the fact that the domain of definition of such a function remains unclarified.  In the particular case of such a system as the one at hand, there is no reason, at present, to expect that any given multiplicity should have a well-defined content.  Since it is far from the case that all properties are heritable, this presents quite a grave problem.  

Recall that we have chosen to build out from the notion of containment---the principle that objects are formed by fixing an interior.  Let an object $a$ have ``content" if there is $b \neq \nll, \neq a$ such that $a \tn b$.  What do we find ``inside" the set of objects $x$ such that $x$ has content?  If the set can be formed, the transitivity of containment becomes destabilizing should there exist objects without content.  Is such a predicate formally unstatable?  It would appear not.  In all nontrivial cases some objects indeed have content; if they all have content, the set is fine.  With that allowance, however, we seem pushed into the boundless possibilities of a hopelessly ill-founded world.  One way to deal with all of this trouble is to suppose that set-theoretical abstraction doesn't care about objects with content in this sense.  If it cares only about objects that lack content, we begin to arrive at the notion of a discrete atomic point as a specialized multiplicity, and a theory of classes based upon it.  

This approach, a retreat into the safety of wellfoundedness, bears a certain resonance, even if one were hypothetically prepared to dismiss the general tendency of modern set theoretical systems as an historic trend.  It remains today an obvious and unavoidable fact (one that no one has ever been quite sure what to make of) that the world of mathematical models is very different, qualitatively speaking, from the world of ``real" objects.  The many-pronged efforts during the early twentieth century to clarify the foundations of mathematics (Hilbert, whose school was perhaps preeminent among these, constantly emphasized the need to divest symbols of their ``meaning" in order to establish perfect rigor) has since been met with the runaway success of modern set theory, logic, and computer science, while the outcome of the movement to clear away the ambiguities of language (led by the logical positivists, significists, and others) remains unclear, having been frustrated by a wave of objections which crested during the postwar period and from which, it seems, it is not likely to recover.  This suggests that the venture to clarify ordinary language and the venture to clarify mathematics are (at least in some significant respects) distinguishable.  Instead of confronting or attempting to repair this dichotomy---and by more tenuous extension, that between mathematical models and empirically observed systems---we would, by enhancing the theoretical significance of the above-mentioned distinction, accept it into our system.  In the endeavor to characterize things which are ``one" (robustly ``single" objects, or \emph{individuals}) contradictions that have beleaguered the study of classes---the paradoxes named after Russell, Curry, Berry, Mirimanoff, Grelling, and so on---shall serve as so many tests to objectively measure candidate definitions.  Should they be resolved, we shall have evidence that a satisfactory abstract definition of the intuitive notion of a point or individual is in reach.  

Our approach will be based loosely upon the notion of primality in a unique factorization domain.\footnote{This structure is discussed in most algebra textbooks, e.g. \citet[\S 3.3]{Hungerford}.}  

\begin{definition}\label{D:Sing}
We say that an object is \emph{single}, $\Sing(a)$, if $(\vas{a} \neq a) \kai ( (a \tn b) \tn (\vas{b} \lor a = b) ).$  That is, $a$ is single if it is a concrete object with no proper concrete content.  A single object is an \emph{individual}.  
\end{definition}

\begin{prop}\label{P:singlehood}
$a=b \tn \Sing(a)=\Sing(b)$.  $\lnot\Sing(\nll)$.  $\lnot\Sing(\dn)$.
\end{prop}

\begin{prop}\label{T:individu}
Individuals are irreducible and prime with respect to $(\kai\,)$.  In symbols:
\begin{gather}
	\Sing(a \kai b) \tn (\vas{a} \lor \vas{b} \lor a = b),  \\
	\Sing(a) \tn (((b \kai c) \tn a) \tn ((b \tn a) \lor (c \tn a))).
\end{gather}
\end{prop}
\begin{proof*}
Let $\Sing(a \kai b)$.  $(a \kai b) \tn a$, so $a = \nll \lor a = (a \kai b)$.  If $a \neq \nll$, then $a = (a \kai b)$; use $(a \kai b) \tn b$.  Now let $\Sing(a)$ and $b \kai c \tn a$.  Suppose $\lnot(b \tn a)$.  $a \tn a \cap b$.  So $a \cap b = a \lor a \cap b = \nll$.  Since $\lnot(b \tn a)$, $a \cap b = \nll$.  So $(b \kai c) \cap a = (b \cap a) \kai (c \cap a) = c \cap a$.  Since $a \neq \nll$ and $(b \kai c) \tn a$, $(b \kai c) \cap a \neq \nll$, so $c \cap a \neq \nll$, so $c \cap a = a$.  So $c \tn (c \cap a) \tn a$.  
\end{proof*}

\begin{definition}
We say that $a$ is an \emph{element} or \emph{member} of an object $b$, $a \in b$, if $\Sing(a)$ and $b \tn a$.  
\end{definition}

\begin{prop}\label{P:marble}
$a \in b \kai c$ implies $a \in b$ or $a \in c$. 
\end{prop}

\begin{prop}\label{P:dnisnot}
$\dn$ is not an element of any object, and nor is $\nll$. 
\end{prop}

Let the reader note, we do not see in $\dn$ the set-theoretical universe---it is a still higher maximum, but it serves only as a kind of universal wastepaper basket.  We will later encounter an object, the ``universe" $V$; it, too, will fail to be an individual in any nontrivial system.  Here, the system echoes the notion, advanced in systems like NBG and MK, that the universe is a proper class, and therefore not an element of any set.

\subsection{Classes}\label{S:3.2}

Because we could ask meaningful questions about classes that have not been constructed---for example, whether (under some formal criteria) they may be constructed or not---we could introduce these objects initially as expressions (say, of the symbolic form $\{x \mid \phi\}$) which dwell in the lattice of objects but have \emph{ab~initio} no containment relationships with any previously defined concrete expression, nor with one another.  Then we might ensure, proceeding abstractly, that the classes so understood behave formally as though they were objects already in some sense defined---possibly by construction.  Where any kind of construction is possible, all conceivable ways of access must converge, and conclude in perfect agreement.  In other words, while the introduction of classes may succeed in generating new descriptions or formulations, it must fail to generate new constructions.  If the class represents a kind of bridge across a sea of complexity, we wish to ensure that it begins and ends {\em somewhere}.  We would not build bridges on a planet with no land, nor if it meant giving up the land for bridges.  

These ideas will resurface several times as we proceed, particularly in \S 4.2. 
Before turning to the formal development, let us consider one more rough analogy.  Consider the sequence $1, 1+1, 1+(1+1), \ldots$ and the sequence $2, 2^{2}, 2^{2^{2}}, \ldots$  The syntax of both expressions is very basic, and one should agree that, whatever these symbols mean relative to one another or to other things, if one of the two sequences could be constructed, then the other could as well.  Most people would also understand the latter expressions to be defined in terms of the former expressions---however, they need not be.  If they were not, then in order to establish equivalence with the ordinary definitions, one would demand that a proof could be supplied showing that we may understand one to be a subsequence of the other, or that a member of the latter sequence takes on a unique location somewhere in the former sequence.  Even if one were not exactly sure \emph{where} in the former sequence, say, the twentieth element of the latter sequence could be correctly placed, one could easily stand about discussing the twentieth element of the latter sequence, more or less as we are doing right now, and proof in hand, refer to it, think of it, regard it, etc., as an element of the former sequence \emph{simpliciter}.  

Now in place of the natural numbers of this example, consider the signs and the formation rules introduced in \S 2. 
We add classes to the system by adding the sign~$\mid$ to the list of basic signs, and to the formation rules of Table \ref{F:frules} the formation rule
\begin{equation*}
	\alpha \, \beta \formation (\alpha\vstop\istop \mid \sbst{\beta}{\alpha\vstop}{\alpha\vstop\istop})
\end{equation*}
We shall axiomatize these new objects in such a way that in every respect they act precisely as though they were in fact members of the original group of well-formed objects---precisely those members that we think of when we imagine the expression factored with respect to $(\kai\,)$.

We define the \emph{class} or \emph{abstract} of variable $x$ and object $\alpha$ (sometimes called the \emph{extension} of $\alpha$ with respect to $x$), denoted $(x\istop \mid \sbst{\alpha}{x}{x\istop})$ or usually more simply with the abbreviation $\clas{x}{\alpha}$, to be a product of $x$, the \emph{class index}, and $\alpha$, the \emph{condition} of the class.  (Though $\alpha$ may range freely over objects, there is little interest in the class if $\alpha$ is not a formula.)  The modified---N.B.---condition $\sbst{\alpha}{x}{x\istop}$ and the modified class index $x\istop$ are subobjects of the class.  We refer to an appearance of an index $x\istop$ as a \emph{bound occurence of $x$}, as usual, or say that $x$ \emph{appears as an index}.  Thus when $x$ simply appears, we sometimes say that $x$ appears \emph{free}.   This allows us in general to do all our thinking in terms of variables, some of which are ``guarded" from substitution processes.

\subsection{Sets}\label{S:B3}

Having discussed the individual and the class, we can introduce the remaining rules and axioms of this section, and begin to develop the rudiments of a theory of sets.  We begin by introducing the simplest class of all:

\begin{definition}
$V = \clas{x}{\nll}$.
\end{definition}

This class may be referred to as the \emph{single-element} or \emph{type-$0$ universe}, or the \emph{universal class}.  The formula $\nll$ which appears as the condition can be thought of as expressing that there are \emph{no conditions} upon members of the class.  

If we select a multiplicity from the object just defined, then we have a set.  

\begin{definition}
An object $a$ is a \emph{set} if $V \tn a$.  If $a$ is a set and $b \tn a$, we say that $a$ is \emph{contained in} $b$, and write $a \lies b$.  If $b$ is also a set, we say that $a$ is a \emph{subset} of $b$.  
\end{definition}

Unless specified otherwise, the letters $\phi$ and $\psi$ shall always denote formulas, and we shall use the variables $A, B, C$, etc. to denote sets.\footnote{As was stated in \S 2, 
all variables are syntactically generic: $\alpha\vstop$ for some object $\alpha$.  We shall not make use of typed variables anywhere in the sequel.}  We define all of the close relatives of $\lies$ and $\in$ in the expected way.\footnote{Because the relation $\lnot(a \in A)$ is rarely used, it might be best to let $a \nin A$ denote that $\Sing(a)$ and $a \not\lies A$.  However, we shall preserve the usual convention regarding slashed relations here.}  Schematically, this gives:
\begin{align}
	\text{\raisebox{-1.1pt}{$\lays$}}\quad 		\ein\quad	\lies\quad	 	&=\quad 		\rise\quad	 	\in\quad		\text{\raisebox{-1.1pt}{$\rais$}}	 	\notag \\
 				\nein\quad	\not\lies\quad 	&\neq\quad 		\not\rise\quad	\nin\quad				 \notag
\end{align}
In expressions these relations will always bind more tightly than the ones introduced in \S 2, 
with the exception of $=$, unless it is denoted by the character $\crn$. 

Passing for a moment to the context of the integers under ordinary multiplication (the model of primality that everyone knows best): we do not regard the integer 1 as a prime, since it is the multiplicative identity, and since it is what mathematicians call a unit.  In a certain sense, it permeates multiplicative space---it is like the solvent in which prime number combinations are dissolved.  Similarly, the basic idea of axiom schema A13 (a version of Cantor's comprehension principle, Frege's rule V, Quine's $\ast 202$, and Zermelo's principle of \emph{Aussonderong}) is to regard $\nll$ as a different kind of object, distinct from ordinary single objects.  It falls outside the domain of things which are targeted for membership in classes---a subset of any set, yet not a \emph{member} of any set.  Granting $\nll$ such special status provides a kind of valve for the release of pressure, which checks the advance of Russell's paradox.

\begin{axiom*}
(Schema) For every well-formed formula $\phi$ and variable $x$,
\begin{equation}\label{E:CompDF}
	a \in \clas{x}{\phi} \crn ( \Sing(a) \kai \phi(a) ).
\end{equation}\vspace{-6pt}
\end{axiom*}

Comprehension schema A13 is presented in a convenient deployable form in (\ref{E:CompDF}), while the formal axiom schema is displayed in Table \ref{F:rax2}. Note, as mentioned above, that since variables, when they are bound, are encased within an index stop, there is no way that a variable appearing in $a$ can be inadvertently bound in deriving $\phi(a)$ from $a$'s membership in $(x\istop \mid \phi(x\istop) )$.  

There is one rule, R3, governing the classes introduced in this section.  It states that under a group of assumptions $\gamma$ which do not depend on a given variable $x$ (i.e., in which the variable $x$ does not appear), if an object $\alpha$ contains an object $\beta$, then the class $\clas{x}{\alpha}$ is contained in the class $\clas{x}{\beta}$.  We point out the immediate derived rule 
\begin{equation}\label{E:dr1}
\begin{array}{c}
	\alpha \tn \beta \\
	\,\clas{x}{\alpha} \lies \clas{x}{\beta}.		\rule{0pt}{10pt}\\
	\raisebox{24.5pt}[0pt][0pt]{$\begin{array}{c}\phantom{xxxxxxxxx.}\\ \hline \end{array}$}
\end{array}\vspace{-6pt}
\end{equation}
The use of $\lies$ here is justified, since for any object $\phi$,
\begin{equation*}
	\clas{x}{\phi} \lies V
\end{equation*}
and R3 itself may be written
\begin{equation}
\begin{array}{c}
	\gamma \tn (\alpha \tn \beta) \\
	\,\gamma \tn \clas{x}{\alpha} \lies \clas{x}{\beta}.		\rule{0pt}{10pt}\\
	\raisebox{24.5pt}[0pt][0pt]{$\begin{array}{c}\phantom{xxxxxxxxxxxxx}\\ \hline \end{array}$}
\end{array}
\end{equation}\vspace{-12pt}

\begin{table}[t]\caption{Postulates of Section 3 (variables $x$, $a$, $b$, objects $\alpha, \beta, \gamma$)}\label{F:rax2}
\begin{center}
\begin{oldtabular}{c}
	\begin{oldtabular}{c}
		R3.  $x$ does not appear in $\gamma$.  \\
		$\begin{array}{c}
			\gamma \tn (\alpha \tn \beta) \\
			\gamma \tn ((x\istop \mid \sbst{\beta}{x}{x\istop}) \tn (x\istop \mid \sbst{\alpha}{x}{x\istop}))	\rule{0pt}{10pt} \\
			\raisebox{24pt}[0pt][0pt]{$\begin{array}{c}\phantom{xxxxxxxxxxxxxxxxxxxxxxxxxx.}\\ \hline \end{array}$}
		\end{array}$\\
	\end{oldtabular}\\
	
	\begin{oldtabular}{l}
		A13 (Comprehension).  objects $\alpha$, variables $x$.  \\
		\quad $a \in (x\istop \mid \sbst{\alpha}{x}{x\istop}) \crn (\Sing(a) \kai \sbst{\alpha}{x}{a})$ \\
		\\
		A14.  variables $x$. \\
		\quad A14a.\;  $a \tn (x\istop \mid x\istop \in a)$  \\
		\quad A14b.\;  $(V \tn a) \tn ((x\istop \mid x\istop \in a) \tn a) $ \\
	\end{oldtabular}
\end{oldtabular}
\end{center}
\end{table}

Since these are the last of our rules, we may now note that the deduction theorem continues to hold, permitting the normal introduction and discharge of assumptions.  We refer to a variable which appears in an open assumption as a \emph{named} variable.  This vocabulary is chosen naturally, because such a variable is thought of as the ``name" of something, and this name is retained as the steps of the proof flow along.  We say that a proof is \emph{mindful of names} if assumptions open during the occurrence of a substitution which has a named variable as target variable appear, upon discharge, with the named variable modified appropriately, and if the rule R3 is used only when the variable to be converted to an index does not appear in any open assumption.  We will assume, from now on, that assumptions are formulas unless stated otherwise.  

\begin{thm}[Deduction theorem]\label{T:deduct2}
Every formal proof with assumptions mindful of names can be converted into a formal proof.  
\end{thm}
\begin{proof*}
We present the verifications to be inserted in Theorem 2.2; 
the setup is identical.  Clearly, the rule R2 can be used to target named variables, provided the proof is mindful of names.  We have to show that instances of R3 can be carried out under the lingering hypothesis $\alpha_i$.  Let an occurence of R3 (using variable $x$) act on step $\alpha_k$, $i < k < i+j$; then $\alpha_k$ is of the form $\gamma \tn (\delta \tn \epsilon)$, where $\gamma, \delta, \epsilon$ are objects, and $x$ does not appear in $\gamma$.  Hence $(\alpha_i \tn \alpha_k) = (\alpha_i \tn (\gamma \tn (\delta \tn \epsilon)))$.  We use this expression to derive a proof (requiring only the postulates of \S 2, 
for which the deduction theorem has been verified) of $(\alpha_i \kai \gamma) \tn (\delta \tn \epsilon)$, apply R3 (noting that the proof is mindful of names), then insert a proof of $(\alpha_i \tn (\gamma \tn \clas{x}{\delta} \lies \clas{x}{\epsilon}))$, allowing the proof to proceed under the assumption $\alpha_i$.  A simple induction argument shows that an arbitrary number of assumptions can be folded up in the same way, to make way for an occurence of R3.\footnote{In fact, because our system is not a resource logic, the open assumptions could simply be treated as an abiding unordered multiplicity.}  This completes the proof.  
\end{proof*}

In certain respects, we expect all sets to behave like deliberately constructed sets.  Hence, one thinks of the set as having a unique decomposition, or of being able to break every set up into its single elements---to ``spread them out" on a universal table top.  Since we have defined sets in terms of a class, $V$, this will be an acceptable notion for us, since all sets will be classifiable, well-founded objects.  At the present stage of development, however, there might exist sets which have subsets, but which do not have any elements.  Suppose, that is, that there exists an object $\ksi$ such that 
\begin{equation*}
\Sing(x) \tn x \nin \ksi,
\end{equation*}
and now suppose that $\ksi$ is a set.  With such an object present, $\clas{x}{x \in (1 \kai 2)}$ could be put equal to $1 \kai 2 \kai \ksi$ without contradiction.  We exclude the possibility of such objects being defined with a version of the principle of extensionality: that every set ``is determined by its elements" \citep{Zermelo1908a}.  The important property is
\begin{equation}
a \tn (b\istop \mid b\istop \in a), 
\end{equation}
which asserts that for any given set $A$, the class $\clas{x}{x \in A}$ is at worst a proper subset of $A$.  For definable sets, the two expressions are equivalent, and any finite union of single objects $a_1\kai \ldots \kai a_n$ is definable (Theorem 3.14). 
For undefinable sets, however, we must enforce a rather imponderable discipline, and thus we include axiom A14b of Table \ref{F:rax2}.

\begin{axiom*}
(Schema) For variables $x$, $\clas{x}{x \in a} \lies a,$ and $a \lies V \tn a \lies \clas{x}{x \in a}$. 
\end{axiom*}\medskip

We can now begin our development in earnest.  

\subsection{Extensionality}\label{S:3.4}

The following theorem is easy to prove in the adopted line of development, and is of central importance to what follows.

\begin{thm}[Extensionality]\label{T:exty}
	$a \lies V \tn ( \clas{x}{x \in a} \lies \clas{x}{x \in b} \tn a \lies b ).$
\end{thm}
\begin{proof*} 
A straightforward application of $\tn$-transitivity, relying on A14.  
\end{proof*}

We can now examine some well-known antinomies.  The Russell set is $s = \clas{x}{x \nin x}$.  A single object $x$ will belong to $s$ if and only if $\lnot(x \in x)$ by A13, but $x \in x$ for all such $x$.  Therefore $s$ contains no individuals.  Therefore (by extensionality) $s$ is empty.\footnote{Though we avoid Russell's paradox, Cantor's theorem on the size of the power set may still be proven (see \S 4.1),
and the theories of computability and decidability appear unaffected, as one might expect.}  A discussion of several similar paradoxically-self-containing sets is found in \citet[pp. 483-4]{Beth1965}.  Extensions of formulas
\begin{equation*}
\call{y_1 \ldots y_n} \quad x \nin y_1 \lor y_1 \nin y_2 \lor \ldots \lor y_{n-1} \nin y_n \lor y_n \nin x
\end{equation*}
(which Quine calls ``epsilon cycles") are under our definitions all equal to $\nll$.  It is not possible under our definitions for any paradox to arise from the relation $s \in s \crn (\Sing(s) \kai \phi(s))$, where $\phi(s)$ implies $s \nin s$.  More generally, we have:

\begin{thm}
$\clas{x}{\dn} = \nll.$
\end{thm}

Curry's paradox arises from considering the extension (with respect to $x$) of the formula $x \in x \tn \phi$, where $\phi$ may be any formula or object.  In naive set theory there is a paradox: let the set be called $c$.  By naive comprehension, $c \in c \crn (c \in c \tn \phi(c))$.  Therefore it is true that $c \in c \tn \phi(c)$.  So $c \in c$.  So $\phi(c)$, which might be any statement at all, is true---all objects, in other words, are trivial.  The A13 comprehension schema also implies that $c \in c \tn \phi(c)$.  However, it does not then follow from A13 that $c \in c$, only that $c \in c$ if $c$ is single.  It then falls to $\phi$ to determine whether in fact $\Sing(c)$ or not.\footnote{In either case, the class is essentially a description which does not affect the underlying ontology.  If $c$ is the sole individual $x$ such that $\phi(x)$---or in other words, if there is but one individual $a$ such that $\phi(a)$---then $c$ is single and $c \in c$, whereas if there are at least two such individuals, or no such individual, it is formally proven that $c$ is not single, and therefore that $c \nin c$.}

The recurring theme of these solutions is that paradoxes are avoided by countenancing multiplicities which are not individuals.  We will see this theme again when we examine other paradoxes in \S 4. 
Next, we verify several basic theorems.  

\begin{thm}\label{T:definablecase}
The class of objects equal to $a$ is $a$ itself if $a$ is single.  The class of all elements of a class is equal to the class itself.  Alphabetic variants are equal; reindexing is permitted.  Briefly:
\begin{gather}
\Sing(a) \tn (\clas{x}{x \in a} = a), \\
\clas{x}{x \in \clas{y}{\phi}} = \clas{y}{\phi}, \\
\clas{x}{\phi} = \clas{y}{\phi(y)}.
\end{gather}
\end{thm}

The proofs of these formulas offer no difficulty.  One simply applies the general method of Theorem 3.12, 
as in, e.g., the following proof:

\begin{thm}\label{T:lor}
For all formulas $\phi$ and $\psi$,
\begin{gather}
\clas{x}{\phi \lor \psi} = \clas{x}{\phi} \kai \clas{x}{\psi},  \label{E:lor-kai} \\
\clas{x}{\phi \kai \psi} = \clas{x}{\phi} \cap \clas{x}{\psi}. \label{E:kai-ka}
\end{gather}
\end{thm}
\begin{proof*}
Let $y \in \clas{x}{\phi \lor \psi}$ (name $y$, and the variables appearing in $\phi$ and $\psi$).  Then by comprehension $\sbst{(\phi \lor \psi)}{x}{y}$, this becomes $\sbst{\phi}{x}{y} \lor \sbst{\psi}{x}{y}$.  Dividing into cases (axiom A6), one concludes $y \in (\clas{x}{\phi} \kai \clas{x}{\psi})$.\footnote{The parentheses required here invite the use of $\cup$ (defined below).}  Discharge the assumption (there are no longer any assumptions, thus no longer any named variables) and apply extensionality.  The other direction is similar.  For (\ref{E:kai-ka}), use A6.
\end{proof*}

\begin{thm}\label{T:cap}
If $A$ and $B$ are sets, $A \cap B = \clas{x}{x \in A \kai x \in B}$.
\end{thm}
\begin{proof*}
As in Theorem 3.15. 
We note that this theorem, like many mathematical theorems, requires a formal stipulation on the involved variables that is often omitted; in this instance, the variable $x$ cannot appear free in the expressions substituted for $A$ and $B$.  For us, however, this restriction is made automatically due to the index-variable distinction drawn at the syntactic level, and the formal statement of the theorem is here given in its entirety.  
\end{proof*}

Since we still lack the ability to form sets of sets, we are still unable to define many familiar devices of set theory, such as the sum set and power set.  However, we may define $A \diff B = \clas{x}{x \in A \kai x \nin B}$, the \emph{difference} of $A$ and $B$, $A \dstnct B = \clas{x}{(x \in A \kai x \nin B)\lor(x \in B \kai x \nin A)}$, the \emph{symmetric difference} or \emph{distinction} of $A$ and $B$, and  $\scomp \! A = \clas{x}{x \nin A}$, the \emph{complement} of $A$.  Under our definitions, the union $A \cup B$ of two sets $A$ and $B$ is once again their federation, the fundamental product $A \kai B$ defined for any two objects $A$ and $B$.  Elementary formulas involving these can now be proven in the usual way, for example: $\scomp\clas{x}{\phi} = \clas{x}{\lnot\phi}$, etc.

\subsection{Quantification}\label{S:3.5}

The principle of Extensionality, first proposed by Frege in 1893, and later by Zermelo in 1908, has been subsequently used to provide a rigorous and intuitively natural basis for theorems throughout mathematics.  We will use it to ground quantification.\footnote{\label{inclusionandabstraction} During the period just prior to the appearance of his work on NF, \citet{QuineIA} experimented with a similar approach, and proved its equivalence to a system of Tarski.  In the system, to the study of which Quine devoted two papers, the simple theory of types is developed using only the connectives of inclusion (containment) and abstraction.  For example, the condional $\phi \tn \psi$ is an abbreviation of the formula $\{x \mid \phi\} \lies \{x \mid \psi\}$, where $x$ does not appear in $\phi$ or in $\psi$.  Since comprehension is one of the basic assumptions in topos theory and related areas, quantification by means of comprehension is commonly used in category-theoretical developments of deductive systems; see, e.g., \citet[p. 70ff]{Bell1988}.  The author is indebted to Andreas Blass for pointing out this and other references in topos theory.}   

\begin{definition}
Given a variable $x$ and formula $\phi$, we write $\call{x}(\phi)$ if $V = \clas{x}{\phi}$, and $\there{x}(\phi)$ if $\nll \neq \clas{x}{\phi}$.\footnote{We will suppress parentheses whenever they are unnecessary.  Our existential quantifier is the double negation of the intuitionistic one.  We might also define symbols for the expressions $\Sing(\clas{x}{\phi})$: there exists a unique individual $x$ such that $\phi$.  Since $\phi(\nll)$ does not imply that $\there{x} \phi$, we might write, say, $\therenll{x} \phi$ to mean $(\there{x} \phi) \lor \phi(\nll)$.  This quantifier is surprisingly rare, however; in mathematical practice, discovery of the existence of structure elements or nonempty sets is the overwhelming norm.  We make note of the presence, in the ordinary language of logic and mathematics, of the words \emph{such that} coupled with the words \emph{there exists} and \emph{the class of}, which suggests---if weakly---a connection between the grammar of the quantifier and the class construction.  But what then, about the incantation ``for all $x$, $\phi(x)$"?  If that hypothesis were true, shouldn't \emph{such that} appear here as well?  Let $\clas{x}{\alpha \mid \phi}$ denote the class $\clas{y}{\there{x}(y \in \alpha \land \phi)}$ (where we insert the needed substitutions should $y$ appear in $\alpha$ or in $\phi$).  Then the universal and existential quantifiers over the scope $\phi$ may be written: \begin{gather*} \nll = \clas{x}{\phi \mid \nll}, \\ \nll \neq \clas{x}{\phi}.\end{gather*}In this form, the scope appears in the condition only in the case of the existential quantifier, prompting \emph{such that $\phi$}, while the condition of the universal quantifier, being empty, is subject to elision; one need not remark the constraint \emph{such that}...well, nothing.  Some might argue that since their scope is confined to $V$, these quantifiers establish a limited or even artificial system of quantification; others, however, might take the view that quantification over objects at the general level has already been defined to satisfaction by the universal rule of substitution R2.}
\end{definition}

Before continuing any further, we need to extend R3 to the case when $x$ may appear in $\gamma$:

\begin{thm}\label{T:R3ex}
$(\clas{x}{\gamma} \lies \clas{x}{\alpha \tn \beta}) \tn ( V \lies \clas{x}{\gamma} \tn \clas{x}{\alpha} \lies \clas{x}{\beta} ).$
\end{thm}
\begin{proof*}
Let $\clas{x}{\gamma} \lies \clas{x}{\alpha \tn \beta}$, $V \lies \clas{x}{\gamma}$.  Then $V = \clas{x}{\alpha \tn \beta}$.  Let $a \in \clas{x}{\alpha}$.  Then $\alpha(a) \tn \beta(a)$, and since $\alpha(a)$, $\beta(a)$ as well, so $a \in \clas{x}{\beta}$.  So $\clas{x}{\alpha} \lies \clas{x}{\beta}$ by extensionality.  
\end{proof*}

Once objects may have classes as subobjects, the management of assumptions becomes something of a puzzle.  Theorem 3.18 
presents a way around this obstacle: it directly follows from it that, given \emph{any} objects $\gamma, \alpha, \beta$, if $\gamma \tn (\alpha \tn \beta)$ has been proven, then one immediately has $(\call{x} \gamma) \tn \clas{x}{\alpha} \lies \clas{x}{\beta}$.  Thus it becomes allowed to convert a named variable to an index, provided the assumptions in which the variable is named are discharged with the introduction of a universal quantifier.  In practice this all falls naturally into place most of the time, but it can incite minor technical fallacies.  

For example, given formulas $\phi$ and $\psi$, the statement $\phi \tn \psi$ should imply that $\clas{x}{\phi} \lies \clas{x}{\psi}$, but as writing $(\phi \tn \psi) \tn \clas{x}{\phi} \lies \clas{x}{\psi}$ produces absurdities by substituting for $x$, this will not be a provable formula unless $x$ does not appear in $\phi$ or in $\psi$.  On the other hand, if for every individual $x$, $\phi \tn \psi$ (a weaker assumption in general), the expected class relation obtains, as has just been shown.  

\citet{Kleene1952} investigates this phenomenon under the rubric \emph{variation} (p. 102ff).  There, one works with a system in which proofs with open assumptions are considered finished formal proofs.  In order to keep his assumptions in order, he stipulates that variables which appear in assumptions are ``held constant" until they are substituted for.  One is free to do so; however, they are then added to a list appended to the proof-theoretical symbol $\vdash$ which indicates that they have been \emph{varied}, and that further substitutions are no longer possible.  His approach may been influenced by the Peanian pasigraphy (by way of the \emph{Principia,} perhaps), in which subscripts of the quantified variables are attached to certain logical connectives in place of quantifier operators.  Having obtained a structural analog to Kleene's approach, we are prepared to proceed.

We next show how to pass from the universal to the particular, and vice versa.  

\begin{thm}[Generalization]\label{T:general}
$\call{a}(\Sing(a) \tn \phi(a)) \tn \call{x}{\phi}.$
\end{thm}
\begin{proof*}
Let $\Sing(a) \tn \phi(a)$.  Suppose $a \in V$; then $\phi(a)$, so $a \in \clas{x}{\phi}$ by comprehension, and therefore $a \in V \tn a \in \clas{x}{\phi}$.  So by extensionality $V \lies \clas{x}{\phi}$.
\end{proof*}

\begin{thm}[Instantiation]\label{T:particular}
$(\Sing(a) \kai \call{x} \phi) \tn \phi(a).$
\end{thm}
\begin{proof*}
1.  Let $\Sing(a) \kai \call{x} \phi.$  2.  $\Sing(a).$  3.  $\Sing(a) \kai \nll.$  4.  $\Sing(a) \kai \sbst{\nll}{x}{a}.$  5.  $a \in V.$  6.  $a \in \clas{x}{\phi}.$  7.  $\phi(a).$
\end{proof*}

In the universal case, we have given a backbone to the quasi-reasoning which echoes in the corridors of mathematics departments, where proofs of universal statements proceed by ``choosing" elements ``arbitrarily", ``freely", sometimes even ``at random".  Existential statements prompt arguments of a similar style: the mathematician argues that he may ``pick one" and insert it into his algebra; if he is subsequently able to derive an expression that does not depend upon which element he picked, he concludes that the statement is true, usually without accounting for the fate of the assumption naming his choice.  This form of reasoning, too, can be filled in with the tools we have defined.  We first pause momentarily to verify the inverse direction:

\begin{thm}[$\exists$-Introduction]\label{T:Eintro}
$a \in \clas{x}{\phi} \tn \there{x} \phi.$
\end{thm}
\begin{proof*}
Argue by contradiction, by concluding from $a = \nll$ that $\Sing(\nll)$.  
\end{proof*}

\begin{thm}[$\exists$-Elimination]\label{T:Eelim}
Let $a$ appear nowhere in formulas $\phi$ and $\psi$, and let $x$ appear nowhere in $\psi$.  Then if, $\call{a}$,
\begin{equation*}
((\there{x} \phi) \kai a \in \clas{x}{\phi}) \tn \psi,
\end{equation*}
then
\begin{equation*}
\there{x} \phi \tn \psi.
\end{equation*}
\end{thm}
We might state the sense of this statement by saying that ``if it has already been proven that $\there{x} \phi \kai \phi(a) \tn \psi$"---here, $a$ is simply an independent variable---``then by using R3 and R1 with the theorem as stated, it is proven that $\there{x} \phi \tn \psi$."\footnote{In Kleene's system, one has $\there{x} \phi, \phi(a) \vdash^{a} \there{x} \phi \tn \psi$.  Notice that $(\there{x} \phi \kai \phi(a) \tn \psi) \tn (\there{x} \phi \tn \psi)$ is false; let $1 \in V$, $a = 1$, $\phi = (x = 1)$, $\psi = (1 \neq 1)$.}

\begin{proof*}
We need the following.
\begin{lem}\label{L:allornothing}
Let $a$ appear nowhere in $\phi$.  Then $(\clas{x}{\phi} \neq \nll) \tn \phi$.
\end{lem}
To prove the lemma, let $\clas{x}{\phi} \neq \nll$, and let $\lnot\phi$.  Suppose (choosing a new variable) that $b \in \clas{x}{\phi}$.  Then $\sbst{\phi}{x}{b}$, that is, $\phi$.  Therefore $\dn$.  Therefore $\sbst{\dn}{x}{b}$, and $b \in \clas{x}{\dn}$.  So $\clas{x}{\phi} \lies \clas{x}{\dn}$ by extensionality.  So $\clas{x}{\phi} = \nll$, contrary to hypothesis.  Discharging gives $\lnot\lnot\phi$, or $\phi$ by the law of excluded middle.  

Let $\phi$ and $\psi$ be defined as stated above.
\begin{stantab}\>\+
	$1.\, \text{Let } ((\there{x} \phi) \kai a \in \clas{x}{\phi}) \tn \psi.$ \\
	$2.\, \text{Let } \there{x} \phi.$ \\
	$3.\, \text{Let } a \in \clas{x}{\phi}.$ \\
	$4.\, \psi$ \\
	$5.\, \sbst{\psi}{x}{a}$ \\
	$6.\, a \in \clas{x}{\psi}$ \\
	$7.\, \text{(3) } a \in \clas{x}{\phi} \tn a \in \clas{x}{\psi}$ \\
	$8.\, \clas{x}{\phi} \lies \clas{x}{\psi}$ (bind $a$) \\
	$9.\, \clas{x}{\psi} = \nll \tn \clas{x}{\phi} = \nll$ \\
	$10.\, (\clas{x}{\psi} = \nll) \tn \dn$ \\
	$11.\, \psi$ (using Lemma 3.23) \\
	$12.\,$ (2), (1) Q.E.D.
\end{stantab}
\end{proof*}

Many readers ought to suspect that a full system of first-order quantification over individuals has been established by these theorems.  We can proceed to develop the identities of classical predicate calculus directly from them, or we can rely instead on the following theorem:  

\begin{thm}\label{T:BQS}
Any theorem provable in classical predicate calculus using the Bernays postulates\footnote{According to \citet[ch. 3]{Hilbert1928}, the postulates defining the universal and existential quantifier in Hilbert-type systems are due to Bernays.} is provable under our postulates \textup{A1-14} and \textup{R1-3}.
\end{thm}
\begin{proof*}
Although the reader may know the four postulates in question well, we should probably review them.  In Kleene's notation, they consist of two axioms: $\call{x} P(x) \tn P(t)$, and $P(t) \tn \there{x} P(x)$, where $t$ is an individual free for $x$ in $P(x)$, along with two rules:
\begin{equation*}
	\begin{array}{c}
		C \tn P(t) \\
		C \tn \call{x} P(x)	\rule{0pt}{10pt}\\
		\raisebox{24.5pt}[0pt][0pt]{$\begin{array}{c}\phantom{xxxxxxxxx}\\ \hline \end{array}$}
	\end{array}
	\qquad\qquad
	\begin{array}{c}
		P(t) \tn C \\
		\there{x} P(x) \tn C	\rule{0pt}{10pt}\\
		\raisebox{24.5pt}[0pt][0pt]{$\begin{array}{c}\phantom{xxxxxxxxx}\\ \hline \end{array}$}
	\end{array}\vspace{-12pt}
\end{equation*}
where the variable $x$ does not appear free in the formula $C$.  

Theorems 3.20 and 3.21 
are analogues of Bernays' axioms, with the proviso that $t$ be an individual converted to one that $t$ be single.\footnote{If we refer to single objects as individuals, then the vocabulary is coincidentally a match.}  The rules may be translated into our system as the formula
\begin{equation*}
	(\call{x}(\phi \tn \psi(x)) \tn (\phi \tn \call{x} \psi)) \kai (\call{x}(\psi(x) \tn \phi) \tn (\there{x}(\psi) \tn \phi))
\end{equation*}
where $\phi$ and $\psi$ are formulas and $x$ does not appear in $\phi$.  This formula may be formally proven using the methods already presented.  
\end{proof*}

\subsection{Elementary Model Theory, Consistency}\label{S:3.6}

We will conclude this section with a short investigation of the metamathematics of the system we have introduced, following \citet{Kleene1952,Simmons}.  Because the new approach to the core syntax of logic adopted in this paper demands some revision of the details of an ordinary consistency proof (while allowing much to stand unchanged), we will proceed in a thorough style, assuming that the reader familiar with logic will easily pick out the relevant changes.  

We begin by introducing revisions to the system that are convenient for the work of this section, since we are only interested in the classical case.  All our previous work goes through essentially unchanged.\footnote{Though in order to undertake this step we must assume primality (Proposition 3.7), 
a principle which as an axiom is unnatural in the author's opinion.  Trying to smooth the wrinkle out causes it to pop up elsewhere; in the balance it seems A5-6 are wanted after all.}
\begin{enumerate}
	\item Cross out axioms A5, A6, and formation rule $\alpha \, \beta \formation (\alpha \cap \beta)$.  
	\item Add abbreviations: $a \lor b$ for $\lnot(\lnot a \kai \lnot b)$, and $a \cap b$ for $\clas{x}{x \in a \kai x \in b}$.
\end{enumerate}

We (momentarily) call an object $\Pi$ a \emph{string} if $\Pi$ may be generated from the following formation rules:
\begin{align*}
	A_i &\formation A_i 						&&\text{ for axiom } A_i \\
	\Pi \, \alpha \, x &\formation \Pi(\alpha/x)	&&\text{ for strings $\Pi$, objects $\alpha$, and variables $x$} \\
	\Pi \, \Sigma &\formation (\Pi \cdot \Sigma) 	&&\text{ for strings $\Pi$ and $\Sigma$} \\
	\Pi \, x &\formation \Pi[x]				&&\text{ for strings $\Pi$ and variables $x$} 
\end{align*}
A string $\Pi$ is subject to the following substitutive one-step reductions:
\begin{align*}
	\alpha(\beta/x) &\vdash \sbst{\alpha}{x}{\beta} 	\\
	(\alpha \cdot (\alpha \tn \beta)) &\vdash \beta		\\
	(\gamma \tn (\alpha \tn \beta))[x] &\vdash (\gamma \tn (\clas{x}{\alpha} \lies \clas{x}{\beta}))
\end{align*}
provided that $x$ does not appear in $\gamma$.  If a string may be reduced to a well-formed object, we call the string a \emph{well-formed proof string} or simply a \emph{proof string}.  We say that an object $\alpha$ is \emph{provable} if there is a well-formed proof string $\Pi$ and a many-step reduction reducing $\Pi$ to $\alpha$.  In this case, we may sometimes call $\alpha$ the \emph{theorem} of $\Pi$, and $\Pi$ a \emph{proof} of $\alpha$, and denote this relation $\Pi \vdash \alpha$.  The \emph{length} of proof string $\Pi$ is defined as usual.

Let $k$ be a nonnegative integer.  A \emph{$k$-structure} $\mathscr{M}_k$ is a metamathematically collected group (hereafter: a set) $(M, P)$, where $M$, the \emph{constants} of $\mathscr{M}_k$, is a set $\nll, 1, 2, \ldots, 2^{k} -1, \dn$.  The \emph{tables} $P$ of $\mathscr{M}_k$ are sets of ordered triplets, often summarized (informally) in the form of written tables.  The tables are a generalization of the notion of a truth table.  It is convenient to correlate the values found in the tables of a $k$-structure to the decomposition of each concrete constant into its base 2 expansion.  For example, in any $k$-structure with $k \geq 4$, $(27,18,\nll)$ will be a member of the $\tn$-table, indicating that $27 \tn 18$ is valid (in the sense defined shortly).  This is because $27 = 2^4 + 2^3 + 2 + 1$, and $18 = 2^4 + 2$ (see Table \ref{F:struct}).  The individuals of each $k$-structure are thus the $k$ constants $1, 2, 4, \ldots, 2^{k-1}$.  We also note that the value of an object $\alpha \tn \beta$ should be, according to $P$, an object having weight 0 (as defined below). 

\emph{Valuation of closed objects.}  Let $\beta$ be a subobject of $\alpha$.  The \emph{class degree} of an instance of $\beta$ with respect to $\alpha$ is the number of classes in $\alpha$ having $\beta$ in its scope, i.e., the number of classes escaped on a walk out of the object beginning at $\beta$.  If an object (formula or concrete) contains no variables, it is said to be \emph{closed}.  The \emph{weight} of a closed object $\alpha$, $\weight(\alpha)$, is defined recursively as follows, letting $\alpha, \beta$ be formal objects: $\weight(\nll) = 0, \weight(\alpha\cstop\,) = \weight(\clas{x}{\alpha}) = \weight(\alpha)$, $\weight(\alpha \kai \beta) = \max(\weight(\alpha) \kai \weight(\beta))$, $\weight(\alpha \tn \beta) = \max(\weight(\alpha) \kai \weight(\beta)) + 1$.  This is a ``metamathematical" function, since its argument is not substitutive.  The procedure we shall adopt for the \emph{evaluation} or \emph{valuation} of a closed object $\alpha$ with respect to a $k$-structure $\mathscr{M}_k$ proceeds in two stages.  

\begin{table}\caption{$k$-Structures}\label{F:struct}
\begin{center}
$\begin{array}{c}
	\begin{array}[t]{c}
		k = 0: \\
		\begin{array}{c|cc}
		\,\;\kai\,\, & \nll & \dn \\
		\nll & \nll & \dn \\
		\dn & \dn & \dn 
		\end{array}\\
		\begin{oldtabular}{c}
		\raisebox{33.4pt}[0pt][0pt]{$\begin{array}{c}\phantom{xxxxxx..}\\ \hline \end{array}$}	
		\end{oldtabular}
		\\
		\begin{array}{c|cc}
		\tn & \nll & \dn \\
		\nll & \nll & \dn \\
		\dn & \nll & \nll 
		\end{array}\\
		\begin{oldtabular}{c}
		\raisebox{33.4pt}[0pt][0pt]{$\begin{array}{c}\phantom{xxxxxx..}\\ \hline \end{array}$}
		\end{oldtabular}
	\end{array} \hspace{-6pt}
	\begin{array}[t]{c}
		k = 2: \\
		\begin{array}{c|ccccc}
		\,\;\kai\,\, & \nll & 1 	& 2 	& 3	 & \dn \\
		\nll 	& \nll & 1 	& 2 	& 3 	& \dn \\
		1 	& 1	& 1	& 3	& 3	& \dn \\
		2	& 2	& 3	& 2	& 3	& \dn \\
		3	& 3	& 3	& 3	& 3	& \dn \\
		\dn	& \dn & \dn & \dn & \dn & \dn 
		\end{array}\\
		\begin{oldtabular}{c}
		\raisebox{66.4pt}[0pt][0pt]{$\begin{array}{c}\phantom{xxxxxxxxxxxxxxxx.}\\ \hline \end{array}$}
		\end{oldtabular}
		\\
		\begin{array}{c|ccccc}
		\tn	& \nll	& 1		& 2		& 3		& \dn \\
		\nll	& \nll	& \dn	& \dn	& \dn	& \dn \\
		1	& \nll	& \nll	& \dn	& \dn	& \dn \\
		2	& \nll	& \dn	& \nll	& \dn	& \dn \\
		3	& \nll	& \nll	& \nll	& \nll	& \dn \\
		\dn	& \nll	& \nll	& \nll	& \nll	& \nll
		\end{array}\\
		\begin{oldtabular}{c}
		\raisebox{66.4pt}[0pt][0pt]{$\begin{array}{c}\phantom{xxxxxxxxxxxxxxxx.}\\ \hline \end{array}$}
		\end{oldtabular}
	\end{array} \hspace{-6pt}
	\begin{array}[t]{c}
		k = 3: \\
		\begin{array}{c|ccccccccc}
		\,\;\kai\,\,	&	\nll	& 1	& 2	& 3	& 4	& 5	& 6	& 7	& \dn \\
		\nll	& 	\nll	& 1	& 2	& 3	& 4	& 5	& 6	& 7	& \dn \\
		1	& 	1	& 1	& 3	& 3	& 5	& 6	& 7	& 7	& \dn \\
		2	& 	2	& 3	& 2	& 3	& 6	& 7	& 6	& 7	& \dn \\
		3	& 	3	& 3	& 3	& 3	& 7	& 7	& 7	& 7	& \dn \\
		4	& 	4	& 5	& 6	& 7	& 4	& 5	& 6	& 7	& \dn \\
		5	& 	5	& 5	& 7	& 7	& 5	& 5	& 7	& 7	& \dn \\
		6	& 	6	& 7	& 6 	& 7	& 6	& 7	& 6	& 7	& \dn \\
		7	& 	7	& 7	& 7 	& 7	& 7	& 7	& 7	& 7	& \dn \\
		\dn	& 	\dn	&\dn&\dn&\dn	&\dn&\dn	&\dn&\dn	&\dn 
		\end{array}\\
		\begin{oldtabular}{c}
		\raisebox{110.4pt}[0pt][0pt]{$\begin{array}{c}\phantom{xxxxxxxxxxxxxxxxxxxxxxxxxxxxx.}\\ \hline \end{array}$}
		\end{oldtabular}\\
		\\
		\text{Etc.}
	\end{array}
\end{array}$
\end{center}
\end{table}

\emph{Direct evaluation}: If every subobject in $\alpha$ has class degree $0$, then the object is class-free.  Class-free objects are evaluated as follows: replace all subobjects having weight 1 with objects having weight 0 by using the tables of $\mathscr{M}_k$.  This shifts the weights of all weighted subobjects by 1.  If necessary, evaluate all multiplicities of weight 0.  Iteration of this procedure eventually yields an object $C$ of weight 0.  

\emph{Declassification}:  If not, then all classes $\pi$ with maximum class degree appearing in $\alpha$ are replaced systematically with a subset of the constants of $M$ which happen to be single, namely $1,2,4,\ldots,2^{k-1}$.  For each $i$, $0 \leq i < k$, $2^i$ is tested in order with the condition $\alpha_\pi$ of $\pi$ (which may be directly evaluated).  If evaluating $\alpha_{\pi}(2^{i})$ yields $\nll$, the constant is included in the multiplicity.  If $\alpha_{\pi}(2^{i})$ takes any other value, then $\nll$ is included in the multiplicity instead.  When this is finished, the multiplicity is itself evaluated, and the class $\pi$ is replaced.  After treating all such classes $\pi$, this step yields an object with strictly smaller class degree.  The entire algorithm is then repeated until $\alpha$ itself is class-free and may be evaluated.  In the case that this procedure yields the product $C$ given $\alpha$, we say that $\alpha$ \emph{takes the value $C$ in $\mathscr{M}_k$}, or that the object is \emph{$C$-valid in $k$}, and write $\smash{\alpha \underset{k}{\Rightarrow} C}$.  If the object is $\nll$-valid, we may say simply that it is \emph{valid}.  

\emph{Objects that are not closed.}  If an object $\alpha$ in which variables appear may be closed by an assignment of its variables to the constants of $M$ (including $\nll$ and $\dn$) such that it takes the value $C$, then we say that the object \emph{can take the value $C$ in $\mathscr{M}_k$}, or that it is \emph{$C$-satisfiable in $k$}.  If the object is $\nll$-satisfiable, we may say simply that it is \emph{satisfiable}.  If under every such assignment $\alpha$ can take no other value except $C$, we again say that $\alpha$ is $C$-valid.  

\emph{Statability.} An object $\alpha$ is \emph{statable in $\mathscr{M}_k$} or \emph{statable in $k$} if no constants except those of $\mathscr{M}_k$ appear in $\alpha$.  $\alpha$ is \emph{provable in $k$}, or $\mathscr{M}_k$-provable, if there is a proof string $\Pi$ that is statable in $k$ (in the obvious extension of the just-stated definition of statability) such that $\Pi \vdash \alpha$.  Note that if an object is provable in $\mathscr{M}_k$, then it is necessarily statable in $\mathscr{M}_k$.  

We refer to $\mathscr{M}_0$ as the \emph{empty structure}.  In accordance with ordinary usage, we refer to class-free formulas statable in the empty structure as \emph{propositions}.  The system restricted to class free formulas (including propositions) statable in the empty structure we call the \emph{extended propositional calculus}.  The system of objects (including arbitrary formulas and concrete objects) statable in the empty structure we call the \emph{extended predicate calculus}.\footnote{Predicates, or predicate letters, can be modeled by variables.  They are functions which take the value $\nll$ or $\dn$ under their inputs, normally taken to be individuals.}

\begin{thm}\label{T:sound}
Any formal object provable in $k$ is valid in $k+l$, for every $l \geq 0$.  In particular, all theorems of the extended predicate calculus are valid in any $k$-structure.  
\end{thm}
\begin{proof*}
We only sketch the proof.  Fix $k$, and let $\Pi$ be statable in $k$.  We wish to show that $\call{l} \geq 0,$ $\smash{\Pi \underset{k+l}{\Rightarrow} \nll}$, that is,
\begin{equation*}
\Pi^{\reduce\,\asg_{k+l}\,\dec_{k+l}\,\ev_{k+l}} = \nll,
\end{equation*}
where the superscripted operators represent (in order) the reduction of $\Pi$ to its theorem (say, $\alpha$), the division of $\alpha$ into a set $\alpha_i$ distinguished by different possible assignments of the variables of $\alpha$ to constants of $\mathscr{M}_{k+l}$, the declassification of the $\alpha_i$ in $k+l$, and finally, the direct evaluation of the declassified $\alpha_i$ using the tables of $\mathscr{M}_{k+l}$.  We proceed by induction on the length $s$ of $\Pi$.  

First, if $s = 0$, then $\Pi = A_i$ is an axiom, statable in the empty structure and thus statable in $k$.  Choosing $l$ arbitrarily, we can verify that in this case $\Pi$ is valid in $l$.  For example, to verify A1, $((c \tn a) \kai (c \tn (a \tn b))) \tn (c \tn b)$, note that if the $c$ is assigned to $\dn$ the resulting formula is valid.  Suppose $c = 2^{c_1}+\ldots+2^{c_n}$, with $0 \leq c_i < l$, and check that in all possible cases, the formula is valid.  This type of inference can be applied to produce verifications of axioms 1, 2, 3, 4, 7, 8, 9, 10, 11, 12, 13, and 14.  We omit the remaining details of this step.  

Next, assume the claim holds for every proof of length $< s$, and let $\Pi$ statable in $k$ have length $s$.  There are three cases to consider.  Let $l \geq 0$.

Case 1.  Let $\Pi = \Sigma[x]$.  First, assume that the theorem of $\Pi$, $\Pi^{\reduce}$, is closed with the possible exception of the variable $x$. We know that the theorem of $\Sigma$, $\Sigma^{\reduce}$, is of the form $\gamma \tn (\delta \tn \epsilon)$, for objects $\gamma, \delta, \epsilon$, and we know that $(\gamma \tn (\delta \tn \epsilon))$ is valid in $k+l$ (hereafter: $(\gamma \tn (\delta \tn \epsilon)) \Rightarrow \nll$).  Since $(\delta \tn \epsilon) \Rightarrow \nll$ (for we may assume that $\gamma \Rightarrow \nll$, since otherwise the valuation of $\Pi^{\reduce}$ is immediately $\nll$), $\delta(2^{c_i}) \Rightarrow \nll$ implies that $\epsilon(2^{c_i}) \Rightarrow \nll$ as well, for any individual $2^{c_i}$ in $\mathscr{M}_{k+l}$.  So if $\clas{x}{\delta} \Rightarrow 2^{c_1} + \ldots + 2^{c_n}$ (as in the previous step), then for every $i$, $\epsilon(2^{c_i}) \Rightarrow \nll$, so $\clas{x}{\delta} \lies \clas{x}{\epsilon} \Rightarrow \nll$, showing that $\Pi$ is valid.  If $\Pi^{\reduce}$ is not closed, the same argument can be repeated on each variable assignment.

Case 2.  Let $\Pi = \Sigma(\beta/x)$.  The validity of $\Sigma^{\reduce}$ implies that assigning any value to the variable $x$, if it appears in $\Sigma^{\reduce}$, the resulting object is valid.  So for any possible value $C$ of $\beta$, $\Sigma^{\reduce}[x/C]$ will clearly also be valid.  Since declassification and valuation both proceed from the leaves to the root of the syntax tree, this proves that $\Pi^{\reduce}$ is valid also in this case.

Case 3.  Let $\Pi = \Sigma_1 \cdot \Sigma_2$.  If $\Sigma_1^{\reduce} = \alpha$, then $\Sigma_2^{\reduce}$ is of the form $\alpha \tn \beta$, for object $\beta$.  Since $\alpha \Rightarrow \nll$, $\beta$ must necessarily also be valid; otherwise $\Sigma_2^{\reduce}$ will not be valid, contrary to hypothesis.  
\end{proof*}

\begin{cor}
No proof string statable in $k$ proves any constant of $\mathscr{M}_k$ except $\nll$.  In particular, $\dn$ is not provable in $k$, for every $k$.
\end{cor}

In the author's opinion, our approach to the distinction between the syntax and semantics of model theory---treating it as a distinction in method instead of a distinction in kind---enhances the model theoretic approach in many respects, simplifying proofs and allowing ideas to more easily commingle with those areas of mathematics where a syntactic/semantic distinction is not implemented.  Other routes to consistency, e.g., via the subformula property, cut elimination, the extended Hauptsatz, and the reduction of the Hilbert-type system to the Gentzen-type systems to which these theorems apply, requires a great deal of technical verification to carry out.  Consistency is more easily provable in natural deduction systems, by using the inversion principle and normalization lemma of \citet{Prawitz1965}, but more effort is required to relate these systems to Hilbert-type systems.   

\section{Second and Higher Order Logic}\label{S:4}

The quantifiers that were introduced in \S 3.5 
are of the first order, insofar as their range is constrained to single objects, or in essence those bodies we usually term individuals.  As yet there is no way, given only the tools defined so far, to form a collection of subsets distinct from their set-theoretic sum, nor is there a way to quantify over such a collection.  The first-order quantifiers of ZF, in stark contrast, range over any set in the cumulative hierarchy.  At this stage, further development requires that we select one of two very broad options.  First, we can bestow upon sets a route to obtaining the property of singlehood, thereby filling the universe $V$ in the way normally done in ZF and similar theories.  Second, following the paradigm of putting eggs into cartons, cartons onto racks, racks onto pallets, pallets into containers, and so on, we can build, atop of $V$, a larger universe of collections.  We present versions of these two approaches in this section, beginning with the latter, and discuss Cantor's theorem and the Burali-Forti paradox.  

\subsection{Genera}\label{S:4.1}

In certain respects the approach of Russell, Ramsey, Quine and other typed systems is quite similar to what follows in this subsection---so much so that the present system can be considered a descendent of theirs.  

We carry out a construction of the \emph{collection of sets}, which we style the \emph{genus}.  Let a pair of signs $\ett$ and $\midI$ be introduced, together with the formation rules
\begin{align*}
	\alpha \, \beta &\formation (\alpha \ett \beta) \\
	\alpha \, \beta &\formation (\alpha\vstop\istop \midI \sbst{\beta}{\alpha\vstop}{\alpha\vstop\istop} ) 
\end{align*}
Thus ($\ett$) is a binary operation on objects, called \emph{type-1 federation}.  Objects $a$ and $b$ are subobjects of the object $a \ett b$, as usual.  In our axiomatization of ($\ett$), we shall rely upon equalities between objects to formally define a product that is identical to union/conjunction/federation ($\kai\,$), except that it does not play the role that federation plays in logic.

\begin{axiom*}  Type-1 federation is associative, commutative, idempotent, substitutive, and has zero element $\nll$ (axioms B1-6 of Table \ref{F:rax3}).  
\end{axiom*}

\begin{definition}
Let $a$ and $b$ be objects.  We write $a \liesI b$, or simply $a \lies b$ if context allows, if $b = a \ett b.$  The relation $\liesI$ is called \emph{type-1 containment} or \emph{first containment}.  When $a \liesI b$, we say that $a$ is \emph{contained in} $b$, or \emph{type-1 contained in}, or \emph{first-contained in} $b$.  
\end{definition}

We can eventually define type-$n$ containment where $n$ is any integer or even any ordinal.  

\begin{thm}
The type-1 containment relation is reflexive, transitive, and antisymmetric.  In addition,
\begin{gather}
	(a \liesI b \kai b \liesI a) \tn a = b. \\
	\nll \liesI a \liesI a \ett b. \\
	(a \liesI b) \tn (a \liesI b \ett c). \\
	(a \liesI c \,\kai\, b \liesI c) \tn (a \ett b \liesI c).
\end{gather}
\end{thm}
\begin{proof*}
To prove that $\liesI$ is transitive, let $a \liesI b$ and $b \liesI c$.  Then $b = a \ett b$ and $c = b \ett c$.  Therefore
\begin{equation*}
c = (a \ett b) \ett c = a \ett (b \ett c) = a \ett c.
\end{equation*}
\end{proof*}

\begin{table}[t]\caption{Postulates of Section 4 (variables $a, b, c$, objects $\alpha, \beta, \gamma$)}\label{F:rax3}
\begin{center}
\begin{oldtabular}{c}
\begin{oldtabular}[t]{lc}
	B1. &  $(a \ett b) \ett c = a \ett (b \ett c)$ \\
	B2. &  $a \ett \nll = a$ \\
	B3. &  $a \ett \dn = \dn$ \\
	B4. &  $a \ett b = b \ett a$ \\
	B5. &  $a \ett a = a$ \\
	B6. &  $(a = b) \tn (a \ett c = b \ett c)$ \\
	B7. & $\Sing(a) \tn \SingI(a)$ \\
	B8. & $\SingI(a) \tn a \lies V$ \\
	B9. & $(\lnot\SingI(a) \land b \neq \nll \land b \neq a) \tn$ \\
	       & $(a \kai b = \dn \land \lnot(a \tn b) \land \lnot(b \tn a))$ \\
	B10. & $a \inI b \ett c \tn (a \inI b) \lor (a \inI c).$ \\
	B11. & $a \inI (x\istop \midI \sbst{\phi}{x}{x\istop}) \crn \SingI(a) \kai \sbst{\phi}{x}{a}$ \\
	\qquad & (objects $\phi$, variables $x$) \\
	B12a. &  $(x\istop \midI x\istop \inI a) \liesI a$  \\
	B12b. & $(a \liesI V_1) \tn (a \liesI (x\istop \midI x\istop \inI a)) $ \\
	\qquad & (variables $x$) \\
\end{oldtabular}
\begin{oldtabular}[t]{c}
\\
R4. $x$ does not appear in $\gamma$ \\
$\begin{array}{c}
	\gamma \tn (\alpha \tn \beta) \\
	\gamma \tn (x\istop \midI \sbst{\alpha}{x}{x\istop}) \liesI (x\istop \midI \sbst{\beta}{x}{x\istop})	\rule{0pt}{10pt}\\
	\raisebox{24pt}[0pt][0pt]{$\begin{array}{c}\phantom{xxxxxxxxxxxxxxxxxxxxxxxxxxx}\\ \hline \end{array}$}
\end{array}$
\\
\\
\\
\begin{oldtabular}{lc}
C1. &  $\Sing(a\shl)$	\\
C2. & $a\shl\ashl = a$	 \\
C3. & \;\,$a=b \tn (a\shl = b\shl \land a\ashl = b\ashl)$ 
\end{oldtabular}
\end{oldtabular}
\end{oldtabular}
\end{center}
\end{table}

\begin{definition}\label{D:SingI}
An object is \emph{type-1 single} or \emph{first-single}, or a \emph{generic individual}, $\SingI(a)$, if it has no generic content.  That is, $\SingI(a) \, \crn \, (b \liesI a \tn (\vas{b} \lor b = a)).$  We write $a \in b$, or $a \inI b$, if $a \liesI b$ and $\SingI(a)$.  A genus that is not first-single is a \emph{proper genus}.  A nonsingle, first-single object is a \emph{plurality}.
\end{definition}

Axioms B7-10 further characterize generic individuals.  Axiom B7 ensures that no proper genus can be single (this principle is exigent in light of Theorem 4.33 
below).  Axiom B8 introduces the simplifying assumption that there are no generic individuals beyond the sets.  Thus the generic individuals are the sets: all individuals together with all pluralities.  Axiom B9 grants meaning to the expressions $a \kai b$ and $a \tn b$ when either $a$ or $b$ is a proper genus.  Axiom B10, perhaps the most important of this group, asserts the primality of first-single objects with respect to $(\ett)$.  

\begin{prop}
$\SingI(a \ett b) \tn (\vas{a} \lor \vas{b} \lor a = b).$
\end{prop}
\begin{proof*}
As in \S 3.1. 
\end{proof*}

\begin{prop}
$a \lies V \tn \SingI(a).$
\end{prop}
\begin{proof*}
Suppose that $\lnot\SingI(a)$.  Then $a$ is a proper genus by Definition 4.29. 
Axiom B9 therefore implies that if $a \lies V$, then $a = (a \kai V)$, and hence $a = \dn$.
\end{proof*}

The rest of the apparatus is quite analogous to previous definitions and postulates.  We define, for every variable $X$ and proposition $\phi$, the \emph{type-1 class of $X$ and $\phi$} or the \emph{class of sets $X$ such that $\phi$}, formally as $(X\istop \midI \sbst{\phi}{X}{X\istop}\,)$.  We shall usually denote this $\clas{X}{\phi}$, or $\clasn{1}{X}{\phi}$, if it is necessary to distinguish the type.  The modified condition $\sbst{\phi}{X}{X\istop}$ and the modified class index $X\istop$ are subobjects of the class $\clas{X}{\phi}$.  These objects are subject to postulates B11, B12, and R4.

\begin{definition}
$V_1 = \clas{X}{\nll}$, the class of sets with no condition.  If the object $a \liesI V_1$, $a$ is said to be a \emph{genus}.  If $b$ is a genus and $a \liesI b$, we say that $a$ is a \emph{subgenus} of $b$.  
\end{definition}

$V_1$ is the \emph{type-1 universe}, or the \emph{universal genus}.  Script letters $\mathcal{A}, \mathcal{B}, \mathcal{C},$ etc. are used to denote genera.  By letting $\calln{1}{X} (\phi) = (V_1 = \clas{X}{\phi})$, $\theren{1}{X} (\phi) = (\nll \neq \clas{X}{\phi})$, we define type-1 quantifiers, which range over sets as well as individuals.  We usually denote these $\call{X}\phi$ and $\there{X}\phi$.  

The curious reader may verify that letting 
\begin{equation*}
	a \capI b = \clas{X}{X \inI a \kai X \inI b}
\end{equation*}
everything remains the same as in \S 2 and \S 3 
at the generic level, and flows from the principle of extensionality much as before.  For example, to verify A6, let $c \liesI a$ and $c \liesI b$.  Then $a \capI b = \clas{X}{X \inI a \ett c \kai X \inI b \ett c} =$ (by primality) $\clas{X}{X \inI c \lor (X \inI a \kai X \inI b)} = c \ett (a \capI b).$  Here, as elsewhere, much depends upon the formula $\clas{X}{\phi \lor \psi} = \clas{X}{\phi} \ett \clas{X}{\psi}$.

We will not do much with genera here (though see \S 5), 
but we can define a few familiar objects.  For every genus $\mathcal{A}$, the \emph{sum set} of $\mathcal{A}$ is $\mathfrak{S}\mspace{-2mu}\mathcal{A} = \clas{x}{\there{X} (x \in X \kai X \in \mathcal{A})}.$  We can also define, for all sets $A$, its \emph{set power} or \emph{subset genus}, the genus of all subsets of $A$: $\mathscr{P} \!A = \clas{X}{X \lies A}.$  Thus we have, for single $a, b, c, \ldots$
\begin{gather}
\mathscr{P}(\nll) = \nll, \notag\\
\mathscr{P}(a) = a, \notag\\
\mathscr{P}(a \kai b) = a \ett b \ett (a \kai b), \notag\\
\vdots \notag
\end{gather}
In general, the power of a finite set of $n$ individuals has $2^{n}-1$ elements.  If this is something less than the reader expected, he or she will note that the cardinality of infinite sets is unaffected by the missing null (as noted above, $\nll$ is not an element of any set); in fact, the power is itself of less theoretical importance without the motivation of Zermelo's principle of \emph{Aussonderong} and the Power Set axiom.  The usual value $2^n$ can be regained with a simple adjustment of the definitions above: instead of having as the type-1 identity $\nll$, one introduces a new identity which might be denoted $\nll_1$, or even $\{\nll\}$.  However, since these constants quickly proliferate, the user of the system is then forced to define scores of symbols all intuitively signifying a trivial structure, a ``set of nothing".  The author's approach is therefore to pursue the simpler of the two options, and to stand prepared should a clear and distinct reason for introducing typed identities arise.  This may well occur, even though it has yet to so far; we should recall the nearly two millenia that passed in the West before 0 was finally distinguished from $\nll$.\footnote{The infinite reduplication of each logically definable class at each type in the simple theory of types (which, as noted in \citet[p. 137]{QuineTT}, is ``particularly strange in the case of the null class") was foremost among the concerns that gave rise to NF.}

We assume (without formal development) the notion of a relation $R$ as a set of ordered pairs, that is, of substitutive products $(a,b)$, for objects $a$ and $b$.  In place of $\clas{y}{(x,y) \lies R}$, we write $R(a)$.  We call a relation $f$ a \emph{function} if $\call{a} \Sing(f(a))$, etc.  

\begin{thm}[Cantor]\label{T:cantorpower}
Let the set $A$ be a plurality.  Then there is no injective function from $A$ onto $\mathscr{P} \!A$.  
\end{thm} 
\begin{proof*} 
Let $A$ be a plurality, let $f$ be an injective mapping of $A$ onto $\mathscr{P} \!A$, and consider the set $Y = \clas{x}{x \nin f(x)}$.  On pain of contradiction, $Y$ must be empty.  In that case (encountered when the Russell set was examined above), an element of $A$ can be mapped to no other value but itself.  But since $A$ is a plurality, this contradicts the hypothesis that $f$ is onto.
\end{proof*}

Because Cantor's theorem is unavoidable, genera cannot be members of $V$.  If genera were allowed to be single, Cantor's theorem would lead directly to Cantor's paradox, because the power of $V$, a subgenus of $V_1$, would be mapped in its entirety into $V$ itself, leading to contradictory results on the cardinality of $V$.  This means, effectively, that it is not possible to form a set of proper genera.  Per axiom B9, such sets are left undefined. 

\subsection{Totalities}\label{S:4.2}

The reader is no doubt familiar with the construction introduced by von Neumann in 1923 defining sets that may be conveniently identified with order types.  This construction is carried out within our system in this section.  The original intent and style of thought is preserved: rather than sets (as given above) we formally generate a class of individuals sharing characteristic properties (which we label \emph{totalities}) and apply abstraction to study them as an unrestricted whole.  We begin by introducing two new symbols, $\shl$ and $\ashl$, with formation rules 
\begin{gather*}
\alpha \formation (\alpha)\shl 		\\ 
\alpha \formation ((\alpha)\shl)\ashl 
\end{gather*}
(we will suppress these parentheses when they are unnecessary).  With these come axioms C1-3:\footnote{As noted above, the product $\land$, equivalent to $\kai$ except for a lower place in the order of operations, may be used at times to suppress unwanted parentheses.  This \label{OoOp} gives rise to the order of operations as follows (where, for added context, ring operations are added): \begin{equation*} \cdot, \kai, +, \cap, \cup, \{=, \lies, \rise\}, \land, \lor, \{ \crn, \tn, \lrn \} \end{equation*}}
\begin{gather*}
\Sing(a\shl)					\\
a\shl\ashl = a					\\
a=b \tn (a\shl = b\shl \land a\ashl = b\ashl)
\end{gather*}
The operation $(\cdot)\shl$ is analogue to the braces of normal set notation.  Following D.~R.~Finkelstein, we call it the \emph{bracing} operation.  The \emph{antibracing} operation $\ashl$ strips away the shell introduced by bracing.  We call objects $a$ such that $a = b\shl$ for some $b$ \emph{totalities}.  Let the ordered pair $(a,b)$ of two objects be the totality $(a\shl \kai (a \kai b)\shl)\shl$; relations may now be defined as usual (as in, for example, \citet{Suppes1960}, or \citet{Kunen1980}).  In particular, let $a < b$ if $a \in b\ashl$.  We say that a set $A$ is \emph{transitive} if $\call{b}, b \in A \tn b\ashl \lies A$, as usual; however, a totality $a$ is \emph{subtransitive} if $\call{b}, b \in a\ashl \tn b\ashl \lies a\ashl$.  

\begin{definition}
An object $a$ is an \emph{ordinal} if it is a subtransitive totality and $a\ashl$ is well-ordered by $<$.  
\end{definition}

The first few ordinals are then $\nll\shl, \nll\shl\shl, (\nll\shl \kai \nll\shl\shl)\shl,$ etc., or (introducing abbreviations) $0, 1, 2, (0 \kai 1 \kai 2)\shl = 3,$ etc.

Forming a totality out of the set of all ordinals leads to the oldest modern set-theoretical paradox.

\begin{thm}\label{T:buraliforti}
Let $\Omega = \clas{x}{x \textup{ is an ordinal}}.$ Then if $\Omega\shl \in V$, then $\nll = \dn$.  
\end{thm}
\begin{proof*}
$\Omega$, the set of all ordinals, is transitive and well-ordered by $\in$.  Therefore $\Omega\shl$ is itself an ordinal, that is, $\Omega\shl \in \Omega$.  But for any $x \in \Omega$, $\lnot(x < x)$, that is, $x \nin x\ashl$.  Therefore $\Omega\shl \nin \Omega$, a contradiction.
\end{proof*}

This is certainly not the first time that this paradoxical result has come to light in this fashion.  It was found by Cantor while he was developing the theory of the transfinite, and it was formally derived by Burali-Forti in a paper which appeared around the turn of the century.  Many years later, it impinged upon Quine when he extended his system NF with the class-existence axioms of \citet{QuineMLrevised}.  After the first edition of the book appeared in 1940, J.B. Rosser discovered that the system it contained (known as ML) was inconsistent.  Quine had posited that all classes having stratified conditions and whose quantifiers could range freely over classes could be considered elements, and made members of his universe $V$.  Having worked out the difficulties related to Cantor's paradox of the largest cardinal, he had been convinced, wrongly as it turned out, that the system would not be troubled by the paradox of Burali-Forti.\footnote{A repair of ML due to Hao Wang was used for the 1951 edition cited in the bibliography.}

After what happened to no less creative and exact a logician than Quine, we must count ourselves fortunate that we noticed the paradox before hastening to press ourselves.  This may be due in some part to the simplicity of the derivation here compared to the rather lengthy derivation required in ML.  The problem with the naive procedure seems to be that the formation rule
\begin{equation*}
	\alpha \formation \alpha\shl
\end{equation*}
and the axiom
\begin{equation*}
	\Sing(a\shl)
\end{equation*}
when combined give rise to an ungovernable feedback when operating in tandem with the class structure introduced in \S 3. 
A class can be chosen and traded in for a new individual; the creation of an individual sparks a revision of the classes.  Given such a potential source of instability, it shouldn't be surprising that a paradox was discovered lurking.  

Let us see if we can distill the logic of Theorem 4.35. 
Since the crucial property is $\lnot(x < x)$, we could consider the class
\begin{equation*}
	\clas{x}{x \nin x\ashl}
\end{equation*}
which contains $\Omega$, and falls to the Russell antinomy.\footnote{Call this set $s$.  $s\shl$ is an object; $s\shl$ is single.  If $s\shl \in s,$ then $s\shl \nin s\shl\ashl,$ i.e., $s\shl \nin s$.  And vice versa: if $s\shl \nin s, s\shl \nin s\shl\ashl,$ i.e., $s\shl \in s$.}  It seems that classes can now be definitions of the non-predicative kind rejected by Poincar\'e.  Having betrayed our intention (of \S 3.2) 
never to use any class as a vehicle for the formation of \emph{new} objects in addition to those already in place before the class has been formed, we have found ourselves again bushwhacked by Russell's paradox.  

The problem may be as old as set theory itself, but we have new notions to employ---the distinct concepts of multiplicity and totality---which allow us hope that the paradox might yet be dispelled.  It must be noted, however, that these notions were introduced already for this purpose by Cantor.  First to arrive at the difficulty, he concluded that $\Omega$ is a multiplicity for which there is no corresponding totality.  He wrote to Dedekind in 1899:  ``If we start from the notion of a definite multiplicity (a system, a totality) of things, it is necessary ... to distinguish two kinds of multiplicities...  For a multiplicity can be such that the assumption that \emph{all} of its elements `are together' leads to a contradiction, so that it is impossible to conceive of the multiplicity as a unity, as `one finished thing'."\footnote{See \citet[pp. 113-7]{Heij}.}  He called such multiplicities \emph{absolutely infinite} or \emph{inconsistent multiplicities}.  

Cantor's interpretation of the difficulty, that the the combustible totalities cannot be formed, or that $\lnot\Sing(\Omega\shl)$, $\lnot\Sing(s\shl)$, etc., is motivated by the intuitive view (held by Cantor) that this construction and all such constructions stand on the cusp of the natural boundaries of the classifiable.  For $\Omega$, in particular, this account is suggestive, given the commonplace associations (from a range of sources upholding varying standards of rigor) between the incomplete infinite (the Cantorian ``absolute" infinite) and the inconceivable or impossible.  In one form or another, this vision of a great divide in the ontology of set theory is shared by almost all of Cantor's descendents, but our work has at least one profound contrast with others: it envisions the boundary brought down, so to speak, from a great height, and made into something more familiar: the distinction between one and many.  

By taking what we have here denoted $<$ instead of what we have here denoted $\in$ as set-theoretical membership, the individual is to this system essentially what the set is to NBG and MK.  Since set theorists have a long experience with the use of proper classes in these systems (whose relationship with ZF has long been well understood), accepting the existence of proper, ``unfinishable", or ``untotalitizable" classes as Cantor recommended will likely seal the present system from further duress.\footnote{For instance, Mirimanoff's paradox is thus resolved.  Mirimanoff considers the set $W$ of all well-founded sets---or rather, the totality $W\shl$ of well-founded \emph{totalities}.  It must be well-founded, since an infinitely descending chain $\ldots < w_n < \ldots w_1 < W\shl$ from $W\shl$ implies an infinite descending chain from $w_1$, but $w_1 < W\shl$ is well-founded.  If $W\shl$ is well-founded, however, then $W\shl < W\shl$, contradiction.  In NBG, it is deduced that $W$ is a proper class---or in Cantor's idiom, $W$ is an inconsistent (we would almost prefer to say ``perfect") multiplicity---one which cannot be gathered into a whole.}  As Fraenkel points out, however, in the preface of \citet[\S 7]{Bernays1958}, a ``vast field of uncertainty still remains" between those totalities which trigger an alarm and those which can be constructively imagined.  In most axiomatizations of set theory (including ZFC, NBG, and MK), this field of uncertainty is pushed back with seemingly innocuous reassurances---$\Sing((a \kai b)\shl)$, $\Sing(\mathscr{P}(x)\shl)$, $\Sing(\bigcup(x)\shl)$, $\Sing((f(A)\shl$)---and other axioms impacting the structure of the transfinite.

There is, in any case, a guiding assurance: the paradox of the largest ordinal number has arisen for us in precisely the way it arose for Quine in the 1940s---because computationally tricky totalities (or elements, in the case of ML) have been introduced to an underlying system which, as Quine put it, ``is a completely safe basal logic...to which more daring structures may be added at the constructor's peril" \citep[p. 135]{QuineEN}.  As an independent form, the system of \S 3 
remains a consistent and adequate setting for basic mathematical intuition---the intuition of aggregates of objects, finite in number, invested with location in a common space.  This lends support to the view that the problem of the Burali-Forti paradox should be diagnosed as an attribute of the construction of na\"ively defined totalities, and does not reflect doubt upon the fundamental principles governing the multiplicity and the class outlined above.  

There likely exists a number of modifications of C1-3 (perhaps along with the formation rules for totalities) by which the difficulties we have observed could be compellingly resolved.  Set theorists have already focused a great and impressive effort upon this pursuit.  One might also pursue an orthogonal route, laying axioms C1-3 aside and carrying out further constructions and experiments leading into different mathematical and logical domains, employing the B axioms, or something similar, when need arises for collections of sets and second order quantifiers.  If we use the system of \S 3 
to stage the integers (beginning along the lines of $\Sing(1), \Sing(a) \tn \Sing(a+1),$ etc., or by defining Church numerals or combinatory terms), the real number system, an elementary topos, or an abstract group, things might all proceed as one might hope and expect, and we could soon find ourselves having access to a language quite indiscernible from that of researchers whose work is formally grounded by other means.  However, because practical attainment of these diaphanous connections might pose unexpected challenges (and these may be instructive and fruitful), we await the findings of further investigations.   

\section{Lower Order Logic, Modulation, Mathematics}\label{S:5}

Having completed the task set out in the introduction, we are prepared to conclude.  However, we hesitate, since readers might now still wish to ask: why go to all the trouble of laying set-theoretical and logical foundations when, to be sure, set-theoretical and logical foundations have already been laid?  There are several ways to answer this objection based on our previous work.  There is the need in science for endless exploration and inquiry; there is the simplified form of our basic definitions and the pull of Occam's razor; there is the  manner in which the elements of the system are well-correlated to natural human mathematical intuition, providing unity and explanatory power; there is the convenient reduction in the complexity of metamathematics provided, in particular, by using syntactically defined classes, syntactically defined indices and variables, and the synoptic rule R3; the manner in which the great nineteenth-century set-theoretical paradoxes have been, we believe, freshly illuminated, and perhaps better understood.  In sum, the construction has so far proven to possess  a number of compelling features.  Our purpose in this section is to add a final observation to this reply felt to be worthy of note, by briefly discussing the manners in which our prior definitions offer practical utility to a working mathematician, and thus gain ground towards what we consider an essential motive for research in foundations.  We ask the reader's forgiveness if in this section we do not pause to define all of our terms perfectly; the previous work ought to be a satisfactory guide to those seeking a more formal treatment.  We begin by generalizing axiom group B slightly.  

\begin{axiom*}  For nonzero integers $n$ and $m<n$ and objects $a$, $b$, and $c$, the \emph{type-n federal product}, or \emph{type-n federation}, is defined, subject to the following axioms:
\begin{equation*}
\begin{tabular}{lcl}
	\textnormal{B1.}\qquad\qquad & $(a \kain{n} b) \kain{n} c = a \kain{n} (b \kain{n} c).$ & \qquad\qquad \\
	\textnormal{B2.}\qquad\qquad & $a \kain{n} \nll = a.$  & \qquad\qquad \\
	\textnormal{B3.}\qquad\qquad & $a \kain{n} \dn = \dn$ & \qquad\qquad \\
	\textnormal{B4.}\qquad\qquad & $a \kain{n} b = b \kain{n} a.$  & \qquad\qquad \\
	\textnormal{B5.}\qquad\qquad & $a \kain{n} a = a.$  & \qquad\qquad \\
	\textnormal{B6.}\qquad\qquad & $(a = b) \tn (a \kain{n} c = b \kain{n} c).$ & \qquad\qquad \\
	\textnormal{B7.}\qquad\qquad & $\Singn{n-1}(a) \tn \Singn{n}(a).$ & \qquad\qquad \\
	\textnormal{B8.}\qquad\qquad & $\Singn{n}(a) \tn a \liesn{n-1} V_{n-1}.$  & \qquad\qquad \\
	\textnormal{B9.}\qquad\qquad & $\lnot\Singn{n}(a) \tn ((b \neq \nll \kai b \neq a) \tn (a \kain{m} b = \dn \kai \lnot(a \liesn{m} b) \kai \lnot(b \liesn{m} a))).$  & \qquad\qquad \\
	\textnormal{B10.}\qquad\qquad & $a \inn{n} b \kain{n} c \tn a \inn{n} b \lor a \inn{n} c.$  & \qquad\qquad 
\end{tabular}
\end{equation*}
\end{axiom*}

The analogous generalization of axioms B11, B12, and R4 is straightforward.  Type-$n$ abstraction is denoted $\clasn{n}{x}{\phi}$, or $\clas{x}{\phi}$ if context allows; elements $a$ of the class must be type-$n$ single, $\Singn{n}(a)$.  We write $\inn{n}$, $\liesn{n}$, etc.; type subscripts may often be suppressed if confusion does not arise as a result.  The type-$n$ universe $\clasn{n}{x}{\nll}$ is denoted $V_n$.  We set $a \kain{0} b = a \kai b$, etc.  

Names tend to lend concreteness to abstract ideas, so we shall say that if...then we call the object a...
\begin{equation} \notag
\begin{tabular}{lc}
	$x \lies V_{-1},$ \qquad\qquad& residue. \\
	$x \lies V_{0},$ \qquad\qquad & set. \\ 
	$x \lies V_{1},$ \qquad\qquad & genus. \\ 
	$x \lies V_{2},$ \qquad\qquad & family. \\
	$x \lies V_{3},$ \qquad\qquad & phylum. \\
	$x \lies V_{4},$ \qquad\qquad & regnum. 
\end{tabular}
\end{equation}
This taxonomic hierarchy is of course taken over from the ranks of the eighteenth century Linnaean system.  The terms \emph{order} and \emph{class} have been removed, since they are frequently used for other purposes in mathematics.\footnote{Nor are they alone: the terms \emph{form,} \emph{variety,} \emph{type,} and \emph{domain} all have currency in taxonomy.}  The terms \emph{set} and \emph{residue} are added since there is an admirable correlation between these and existing mathematical notions.  Since the term \emph{family} is often used for a set of collections (as for example in \citet{RudinRCA}), there is another agreement with standard terminology at the $V_{2}$ rank.  The term \emph{species} might be suggested over \emph{set}, since Aristotle's original binomial taxonomy can then be recognized in the generalized system, and  since that term is generally used by Brouwer in his papers on intuitionistic set theory (see, e.g., \citet{Brouwer1923}).  It also suggests itself since it is  useful to have the term \emph{set} available at the metalevel (as we saw in \S 3.5). 
Nevertheless, we shall continue to refer to generic individuals as sets.  

\begin{definition}
Let $A = \clasn{n}{x}{x \inn{n} A} = \clasn{n}{x}{\phi}$.  The \emph{modulation} of $A$ is defined to be the expression
\begin{equation*}
\clasn{n-1}{x}{x \inn{n} A} = \clasn{n-1}{x}{\phi},
\end{equation*}
conveniently denoted {\rm $\mdln{A}$}.
\end{definition}

The inverse operation, or \emph{antimodulation} of $A$ may also be defined in the natural way.  Thus from a given set $A$, one obtains a residue through modulation, and a genus through antimodulation.  Among the two, the modulation (which may be pronounced \emph{modulo} $A$, or simply \emph{mod} $A$) serves particularly well as a basic tool in algebra.  Let $R$ be a set, and let $(+), (\cdot)$ be products defined on $R$, with $a+b$, $ab = a \cdot b \in R$ for $a,b \in R$.  Let
\begin{gather}
	a + (b \kain{n} c) = (a + b) \kain{n} (a + c), \\
	a (b \kain{n} c) = ab \kain{n} ac, 
\end{gather}
so that $(+)$ and $(\cdot)$ are distributive over the entire federal hierarchy.  Let $(R,+,\cdot, 0, 1)$ be a ring with unity, and let $S \lies R$.  The question of when the set $R+\mdln{S}$ forms a ring may be answered by examining the elementary calculations 
\begin{gather}
(a + \mdln{S}) + (b + \mdln{S}) = a + b + 2 \mdln{S},  \label{E:44-1}\\
(a + \mdln{S})(b + \mdln{S}) = ab + a \mdln{S} + b \mdln{S} + \mdln{S}\mdln{S}. \label{E:44-2}
\end{gather}
We wish to obtain a set of conditions under which these expressions become
\begin{gather} 
(a + \mdln{S}) + (b + \mdln{S}) = (a + b) + \mdln{S}, \label{E:4-1}\\ 
(a + \mdln{S})(b + \mdln{S}) = ab + \mdln{S},  \label{E:4-2}
\end{gather}
which correlate multiplication and addition with and without the modulated factor $S$.  (\ref{E:4-1}) will follow from (\ref{E:44-1}) when $S$ is additively closed and contains 0.  A solution to the problem of satisfying the constraint of (\ref{E:4-2}) is to allow that for all $a \in S, aS=S$, and that $S$ be multiplicatively closed, so that $S^2 \lies S$ (remembering $0 \in S$).  Thus if $a \in S$, $2a, 3a, \ldots$ and $a^2, a^3, \ldots \in S$ as well.  Moreover $-a \in S$ since $-a \in (-1)S$, giving an additive group, and the multiplicative closure requirement can similarly be dropped.  We call such a structure an \emph{ideal} inside the ring.  

The reader should note that the ideal concept has been defined neither dogmatically nor informally.  Rather, it has been directly detected and recovered, by means of a brief, real-time examination of algebraic structure.  The style of calculation can be continued to yield elegant direct proofs of the isomorphism theorems, and generally to facilitate investigation within the theory of homomorphisms and beyond.  

Consider the case $R = \Z$.  Each ideal has the form $m\Z = \ldots \kai -2m \kai -m \kai 0 \kai m \kai 2m \kai \ldots \; $ for some $m \in \Z$.  The study of equations in $\Z$ carrying the factor $m\mdln{\,\Z}$ was first undertaken systematically by \citet{Gauss}, who wrote 
\begin{equation*}
	a \equiv b \;\, (\text{mod.} \; m)
\end{equation*}
for the relation
\begin{equation*}
a + m\mdln{\,\Z} = b + m\mdln{\,\Z},
\end{equation*}
or equivalently, 
\begin{equation*}
a \in b + m\Z.
\end{equation*}
Thus the formal expressions we obtain using the algebraic set bear a striking resemblance to the classical notation coined by Gauss.  

We present a proof of a well-known elementary theorem in the calculus of residues, a proof that will in fact be similar to the one given by Gauss (\emph{op. cit.}, pp. 10-13).  Because the symbol $\Z$ recurs frequently, it will be suppressed under modulation, thus $m\mdlo = m\mdln{\,\Z}$.  We will work with equalities and containments between residues, i.e., federal products $a + m \mdln{\,\Z}$.  

Consider the problem of finding the integer solutions $x$ of the equation $ax = b$ modulo an integer $m$, that is, finding the complete solution $x \lies \Z$ of the expression
\begin{equation}\label{Z-1}
ax + m\mdlo = b + m\mdlo
\end{equation}
where $a,b,m \in \Z$.  In order to bring the expression into a shorter and more standard form, we write (\ref{Z-1}) as
\begin{equation}\label{Z-2}
{ax = b}^{\,\ssplus m\ssmdlo}
\end{equation}
to indicate that the superscripted factor is applied on both sides of the equation.  Let $d > 0$ be the greatest common divisor of $a$ and $m$.  Adding\footnote{This operation is permitted; however, ideals do not have additive inverses, so one may not subtract them out.} the residue $d\mdlo$ to both sides of (\ref{Z-2}) gives 
\begin{equation}
{ax = b}^{\,\ssplus m\ssmdlo \ssplus d\ssmdlo}
\end{equation}
but since $ax+d\mdlo = d\mdlo$ and $m\mdlo+d\mdlo = d\mdlo$, this becomes
\begin{equation}
{0 = b}^{\,\ssplus d\ssmdlo}
\end{equation}
or
\begin{equation}
d\mdlo = b + d\mdlo.
\end{equation}
Selecting a member from the right hand side gives
\begin{equation}
d\mdlo \ni b.
\end{equation}
Antimodulation gives $b \in d\Z$, i.e., $d$ divides $b$.  Dividing by $d$ gives
\begin{equation}
{a'x = b'}^{\,\ssplus m'\ssmdlo}
\end{equation}
where $a'$ is relatively prime to $m'$.   Since $a\Z + m\Z = d\Z$ (which we do not pause to prove), $b' \in \Z = a'\Z + m'\Z$, and therefore
\begin{equation}
{b' \in a'\Z}^{\,\ssplus m'\ssmdlo}
\end{equation}
that is,
\begin{equation}\label{Z-3}
{b' \in a'(0 \kai 1 \kai 2 \kai \ldots \kai m'-1)}^{\,\ssplus m'\ssmdlo.}
\end{equation}
By multiplying (\ref{Z-3}) through by $d$, we see that the solution of (\ref{Z-2}) is nonempty.  Since ${a' \neq 0}^{\,\ssplus m'\ssmdlo}$, this solution is single modulo $m'$.\footnote{We do not intend to flout tradition: a solution \emph{exists}, and is \emph{unique} modulo $m'$.  Our intent here is to point out a certain agreement between the shared notions of mathematics and notions we have developed in the preceding sections while working in foundations.}  We may thus write the solution set in the form $x_0 + m'\Z$, for some single solution $x_0$.  It remains only to determine the form of the solutions in the modulus $m'$.  This is easy: we proceed by extracting the coarser modulus: 
\begin{align}
x 	&= x_0 + m'\Z 	\notag\\
	&= x_0 + m'\Z + m\Z \notag\\
	&= x_0 + m'(0 \kai 1 \kai \ldots \kai d-1) + m\Z.
\end{align}
We therefore conclude that, when $d$ divides $b$, there are infinitely many solutions of (\ref{Z-2}) in the set of integers, $d$ solutions modulo $m$, and a single solution modulo $m'$.  The argument has a more elementary character than a standard proof, since it relies upon natural extensions of familiar algebraic operations and has a very basic quantificational structure.  We let this suffice as a simple demonstration of the computational methods organized by instruments developed in the foregoing sections, since we intend to continue this discussion in work now in progress.   

\section{Acknowledgements}
This study began while the author was a student at the University of Georgia advised by Brad Bassler, Robert Varley, David Edwards, and Edward Halper.  The author would like to thank them, as well as David Finkelstein and the anonymous referees, for all their help in editing and improving earlier versions of this paper.  

\bibliographystyle{rsl}
\bibliography{settl_final}

\vspace*{10pt}
\address{DEPARTMENT OF MATHEMATICS\\
\hspace*{9pt}LOUISIANA STATE UNIVERSITY\\
\hspace*{18pt}BATON ROUGE, LA 70803\\
{\it E-mail}: lschoe2@lsu.edu\\
{\it URL}: \url{http://www.math.lsu.edu/~lschoe2}}
\clearpage
\end{document}